\documentclass[11pt,a4paper]{article}
\usepackage[german,english]{babel}
\usepackage{latexsym}
\usepackage{amsmath, amsfonts}
\usepackage{amssymb}	
\usepackage{amsthm}	
\usepackage{ifthen}
\usepackage{hyperref} 
\usepackage{graphicx} 
\usepackage{verbatim} 

\author{Michael Klotz\footnote{Technische Universit\"at Darmstadt, Schlo\ss gartenstra\ss e 7, D-64289 Darmstadt, Deutschland,\newline klotz@mathematik.tu-darmstadt.de}}
\date{}
\title{Banach Symmetric Spaces}

\newenvironment{definition}{\begin{Definition}}{\end{Definition}\par}	
	\theoremstyle{definition} 
	\newtheorem{Definition}{Definition}[section]
\newenvironment{theorem}{\begin{Theorem}}{\end{Theorem}\par}
	\theoremstyle{plain} 
	\newtheorem{Theorem}[Definition]{Theorem}	
\newenvironment{proposition}{\begin{Proposition}}{\end{Proposition}\par}
	\newtheorem{Proposition}[Definition]{Proposition}	
\newenvironment{lemma}{\begin{Lemma}}{\end{Lemma}\par}
	\newtheorem{Lemma}[Definition]{Lemma}
\newenvironment{corollary}{\begin{Corollary}}{\end{Corollary}\par}
	\newtheorem{Corollary}[Definition]{Corollary}	
\newenvironment{example}{\begin{Example}}{\end{Example}\par}
	\theoremstyle{definition}
	\newtheorem{Example}[Definition]{Example}
\newenvironment{remark}{\begin{Remark}}{\end{Remark}\par}
	\theoremstyle{definition}
	\newtheorem{Remark}[Definition]{Remark}

	\theoremstyle{definition}
	\newtheorem{Problem}[Definition]{Problem}
\renewenvironment{proof}{\noindent\textbf {Proof:}}{\hspace*{\fill} $\Box$\par\vspace{0.7ex}}

\newenvironment{acknowledgements}{\section*{Acknowledgements}}{}

\newcommand{\temp}{0} 
\newenvironment{temparraystretch}[1] 
	{\renewcommand{\temp}{\arraystretch} \renewcommand{\arraystretch}{#1}}
	{\renewcommand{\arraystretch}{\temp}}

\newcommand{\NN}{\mathbb N}
\newcommand{\RR}{\mathbb R}
\newcommand{\1}{{\bf 1}}

\newcommand{\calD}{{\cal D}}
\newcommand{\calF}{{\cal F}}
\newcommand{\calL}{{\cal L}}
\newcommand{\calV}{{\cal V}}
\newcommand{\extd}{\mbox{\rm d}} 
\newcommand{\gl}{\mathfrak{gl}}

\DeclareMathOperator{\Ad}{Ad}
\DeclareMathOperator{\Alt}{Alt}
\DeclareMathOperator{\Aut}{Aut}
\DeclareMathOperator{\conv}{conv}
\DeclareMathOperator{\Der}{Der}
\DeclareMathOperator{\Diff}{Diff}
\DeclareMathOperator{\ev}{ev}
\DeclareMathOperator{\flow}{Fl}
\DeclareMathOperator{\Fr}{Fr}
\DeclareMathOperator{\GL}{GL}
\DeclareMathOperator{\Gr}{Gr}
\DeclareMathOperator{\id}{id}
\DeclareMathOperator{\im}{im}
\DeclareMathOperator{\Invol}{Invol}
\DeclareMathOperator{\Iso}{Iso}
\DeclareMathOperator{\Kill}{Kill}
\DeclareMathOperator{\Lts}{Lts}
\DeclareMathOperator{\sgn}{sgn}
\DeclareMathOperator{\Tor}{Tor}

\begin{document}
\maketitle
%
%
%
%
\begin{abstract}
	A Banach symmetric space in the sense of O.~Loos is a smooth Banach manifold $M$ endowed with a multiplication map $\mu\colon M\times M\rightarrow M$ such that each left multiplication map $\mu_x:=\mu(x,\cdot)$ (with $x\in M$) is an involutive automorphism of $(M,\mu)$ with the isolated fixed point $x$. We show that morphisms of Lie triple systems of symmetric spaces can be uniquely integrated provided the first manifold is 1-connected. The problem is attacked by showing that a continuous linear map between tangent spaces of affine Banach manifolds with parallel torsion and curvature is integrable to an affine map if it intertwines the torsion and curvature tensors provided the first manifold is 1-connected and the second one is geodesically complete. Further, we show that the automorphism group of a connected Banach symmetric space $M$ can be turned into a Banach--Lie group acting smoothly and transitively on $M$. In particular, $M$ is a Banach homogeneous space.
	
	
\end{abstract}

%
%
%
\section{Introduction}
\label{introduction}
A symmetric space in the sense of O.~Loos (cf.\ \cite{Loo69}) is a smooth manifold $M$ endowed with a multiplication map $\mu\colon M\times M\rightarrow M$ such that each left multiplication map $\mu_x:=\mu(x,\cdot)$ (with $x\in M$) is an involutive automorphism of $(M,\mu)$ with the isolated fixed point $x$. 

The tangent bundle $(TM,T\mu)$ of a symmetric space $(M,\mu)$ is again a symmetric space so that a smooth vector field on $M$ is called a derivation if it is a morphism of symmetric spaces. Each symmetry $\mu_x$ induces an involutive automorphism of the Lie algebra $\Der(M,\mu)$ of derivations. This provides an additional structure on the tangent space $T_xM$, namely a Lie triple system, and we obtain a functor $\Lts$ from the category of pointed symmetric spaces to the category of Lie triple systems.

The purpose of this paper is to start working towards a Lie theory of symmetric spaces modelled on Banach spaces. From \cite{Nee02}, \cite{Lan01} and \cite{Ber08}, we extract some basic material on infinite-dimensional symmetric spaces.

In the finite-dimensional case, O.~Loos claims that morphisms of Lie triple systems of symmetric spaces can be uniquely integrated provided the first manifold is 1-connected (cf.\ \cite{Loo69}). 
A connected symmetric space carries a natural affine connection such that morphisms coincide with affine maps and such that morphisms of appendant Lie triple systems coincide with curvature intertwining maps. Therefore, the integrability problem can be translated into the language of affine connections. Symmetric spaces are torsionfree and geodesically complete and possess parallel curvature.

Given a 1-connected affine manifold $M_1$ and a geodesically complete affine manifold $M_2$ with base points $b_1$ and $b_2$, respectively, O.~Loos claims that a torsion and curvature intertwining map between the tangent spaces at the base points can be uniquely integrated to an affine map from $M_1$ to $M_2$. Regrettably, his line of argument seems to be incomplete. He actually shows only a local extension.

In \cite{Amb56}, W.~Ambrose gives a theorem about the integration of an isometric isomorphism between tangent spaces of complete 1-connected Riemannian manifolds in the finite-dimensional case. As a sequel to it, in \cite{Hic59}, N.~Hicks deals with the case of affine manifolds. Their work is based on the work of \'E.~Cartan (cf.\ \cite{Car46}) and their results are known as the theorem of Cartan--Ambrose--Hicks. This is also dealt with in \cite{KN63} and in \cite{Wol67} and is generalized in \cite{PR02} and \cite{BH89} to the case of maps that are not necessarily isomorphisms. In \cite{PR02}, the theorem is proved by constructing affine maps via their graph and an existence theorem for affine submanifolds. In \cite{BH89}, the maps are constructed between the frame bundles over the manifolds rather than between the manifolds themselves, but likewise with Cartan's `technique of the graph'.

In this paper, we generalize in a first step the arguments in \cite{Loo69} to the Banach case to show that a continuous linear map $A\colon T_{b_1}M_1\rightarrow T_{b_2}M_2$ between tangent spaces of affine Banach manifolds $M_1$ and $M_2$ with parallel torsion and curvature is locally integrable to an affine map if it intertwines these tensors.
In a second step, we globalize this theorem by extension along piecewise geodesics and obtain that the map $A$ is uniquely integrable to an affine map $f\colon M_1\rightarrow M_2$ provided $M_1$ is 1-connected and $M_2$ is geodesically complete.

Applying these results to Banach symmetric spaces, we see that morphisms of Lie triple systems of symmetric spaces can be uniquely integrated provided the first manifold is 1-connected.
As a further important result we show that the automorphism group of a connected Banach symmetric space $M$ is a Banach--Lie group acting smoothly and transitively on $M$. In particular, $M$ is a Banach homogeneous space. More precisely, we have $M\cong G/G_b$, where the stabilizer $G_b$ for a point $b\in M$ is an open subgroup of the group of fixed points in $G$ for the involution $\sigma$ on $G$ given by $\sigma(g):=\mu_b\circ g\circ \mu_b$ with the symmetry $\mu_b$ at $b$ (cf.\ \cite[Ex.~3.9]{Nee02} for homogeneous symmetric spaces). This is a simple consequence of \cite{Klo09a}, where the author shows that the automorphism group of a connected geodesically complete affine Banach manifold $M$ can be turned into a Banach--Lie group acting smoothly on $M$.
\tableofcontents
%
%
%
%
%
\section{Affine Banach Manifolds}
\label{sec:affineManifolds}
In this section, we first collect a number of definitions and properties concerning affine connections on smooth Banach manifolds. Afterwards, we cite relevant statements concerning affine and infinitesimal affine automorphisms: Given a connected affine Banach manifold $(M,\nabla)$, the Lie algebra $\Kill(M,\nabla)$ of infinitesimal automorphisms can be naturally turned into a Banach--Lie algebra. If $M$ is geodesically complete, then $\Kill(M,\nabla)$ consists of complete vector fields and the automorphism group $\Aut(M,\nabla)$ can be turned into Banach--Lie group acting smoothly on $M$.
\subsection{Affine Connections on the Tangent Bundle}
\label{sec:affineConnections}
Let $M$ be a smooth Banach manifold, $\pi\colon TM\rightarrow M$ be the natural projection of its tangent bundle and $\pi_{TM}\colon TTM \rightarrow TM$ be the natural projection of the tangent bundle of $TM$. Also the map $T\pi\colon TTM\rightarrow TM$ makes $TTM$ a vector bundle over $TM$ (cf.\ \cite[p.~104]{Lan01}). The composition $\pi\circ\pi_{TM}$ turns $TTM$ into a fiber bundle over $M$.

An \emph{affine connection on $TM$} is a morphism $B\colon TM\oplus TM \rightarrow TTM$ of fiber bundles over $M$ such that $(\pi_{TM},T\pi)\circ B=\id_{TM\oplus TM}$ and such that $B$ is bilinear, i.e., for each $x\in M$, $B_x\colon T_xM\oplus T_xM \rightarrow TTM$ is bilinear. Note that $B_x(v,\cdot)\colon T_xM \rightarrow T_v(TM)$ and $B_x(\cdot,w)\colon T_xM \rightarrow (T\pi)^{-1}(w)$ are indeed maps between Banach spaces. The pair $(M,B)$ is called an \emph{affine Banach manifold}.

In a chart $\varphi\colon U\rightarrow V\subseteq E$, an affine connection $B$ can be written as
$$\begin{array}{cccc}
	TT\varphi\circ B\circ (T\varphi\oplus T\varphi)^{-1}\colon & TV\oplus TV=V\times E\times E   & \rightarrow & TTV=V\times E\times E\times E\\
	& (x,v,w) & \mapsto & (x,v,w,B^\varphi_x(v,w))
\end{array}$$
with a smooth map
$B^\varphi\colon V \rightarrow \calL^2(E,E)$
from $V$ into the space of continuous bilinear maps $E\times E \rightarrow E$, which we call a \emph{local representation of $B$}. Considering two charts $\varphi_1$ and $\varphi_2$, the change of variable formula for the transition map $h:=\varphi_2\circ\varphi_1^{-1}$ is given by
$$B^{\varphi_2}_{h(x)}(dh(x)(v),dh(x)(w))=d^2h(x)(v,w)+dh(x)(B^{\varphi_1}_{x}(v,w)).$$

An affine connection can also be given by a \emph{covariant derivative} $\nabla$, i.e., by a collection $(\nabla^U)_{U \subseteq M \, \mbox{\scriptsize open}}$ of $\RR$-bilinear maps
$$\nabla^U\colon {\cal V}(U) \times {\cal V}(U) \rightarrow {\cal V}(U),\ (\xi,\eta) \mapsto (\nabla^U)_\xi\eta $$
satisfying the conditions
\begin{enumerate}
	\item[(1)] $(\nabla^U)_{f\xi}\eta = f(\nabla^U)_\xi\eta$ \quad ($C^\infty(U)$-linearity in the first variable)
	\item[(2)] $(\nabla^U)_\xi(f\eta) = (\xi.f)\eta + f(\nabla^U)_\xi\eta$ \quad (derivation property)
\end{enumerate}
for all $\xi,\eta\in{\cal V}(U)$ and smooth functions $f\in C^\infty(U)$ such that the maps $\nabla^U$ are compatible in the sense that
$$((\nabla^{U_1})_\xi\eta)|_{U_2} = (\nabla^{U_2})_{\xi|_{U_2}}\eta|_{U_2}$$
for all $\xi,\eta\in{\cal V}(U_1)$, $U_2\subseteq U_1\subseteq M$.
In the following, we shall often suppress the index set $U$ by writing $$\nabla\!_\xi\eta:= (\nabla^U)_\xi\eta$$
for all $\xi,\eta\in {\cal V}(U)$.

There is a one-to-one correspondence between affine connections and covariant derivatives. It is determined by the local formula
$$(\nabla\!_\xi\eta)^\varphi(x) = d\eta^\varphi(x)(\xi^\varphi(x)) - B^\varphi_x(\eta^\varphi(x),\xi^\varphi(x)),$$
where $(\nabla\!_\xi\eta)^\varphi$, $\eta^\varphi$ and $\xi^\varphi$ denote the local representations of the vector fields.
As far as the vector field $\xi$ is concerned, $(\nabla\!_\xi\eta)(x)$ only depends on $\xi(x)$. Therefore, it make sense to define $\nabla\!_v\eta$ for vectors $v$.


Given a smooth curve $\alpha\colon J \rightarrow M$, let $\gamma\colon J \rightarrow TM$ be a lift of $\alpha$ to $TM$, i.e., a curve on $TM$ satisfying $\pi\circ\gamma=\alpha$. The {\em derivative of $\gamma$ along $\alpha$} is the unique lift $\nabla\!_{\alpha^\prime}\gamma$ of $\alpha$ to $TM$ that in a chart $\varphi\colon U \rightarrow V\subseteq E$ has the expression
$$(\nabla\!_{\alpha^\prime}\gamma)^\varphi(t) \ =\ (\gamma^\varphi)^\prime(t)-B^\varphi_{\alpha^\varphi(t)}(\gamma^\varphi(t),(\alpha^\varphi)^\prime(t)).$$
We also use the notation $\nabla\!_{\alpha^\prime(t)}\gamma$. A lift $\gamma$ of $\alpha$ is said to be {\em $\alpha$-parallel} if $\nabla\!_{\alpha^\prime}\gamma = 0$.

An affine connection induces {\em parallel transport} along smooth curves. For a curve \linebreak $\alpha\colon J \rightarrow M$ and $t_0,t_1 \in J$, we denote it by
$$P_{t_0}^{t_1}(\alpha)\colon T_{\alpha(t_0)}M \rightarrow T_{\alpha(t_1)}M.$$
It is a topological linear isomorphism and is defined by the property that for each $v\in T_{\alpha(t_0)}M$, the map
$$\gamma_v:=P_{t_0}^{(\cdot)}(\alpha)(v)\colon J \rightarrow TM$$
is the unique curve in $TM$ that is $\alpha$-parallel and satisfies
$\gamma_v(t_0)=v$.
In any chart\linebreak $\varphi\colon U\rightarrow V\subseteq E$, it then satisfies the linear differential equation
$$(\gamma_v^\varphi)^\prime(t) \ =\ B^\varphi_{\alpha^\varphi(t)}(\gamma_v^\varphi(t),(\alpha^\varphi)^\prime(t))$$
and it is uniquely determined by satisfying these equations
for a collection of charts covering the curve $\alpha$ and by satisfying the initial condition $\gamma_v(t_0)=v$.
Along piecewise smooth curves, we can define parallel transport, too, by composing it piecewise.

A {\em geodesic} is a curve $\alpha$ in $M$ whose derivative $\alpha^\prime$ is $\alpha$-parallel, i.e.,
$\nabla\!_{\alpha^\prime}\alpha^\prime = 0$.
For each $v\in T_xM,\, x\in M$, for which the unique maximal geodesic $\alpha_v\colon J \rightarrow TM$ with $\alpha_v^\prime(0)=v$ satisfies $1\in J$, we define
$$\exp(v):=\exp_x(v) := \alpha_v(1).$$
We denote the open domains of $\exp$ and $\exp_x$ by ${\cal D}_{\exp}\subseteq TM$ and ${\cal D}_{\exp_x}\subseteq T_xM$, respectively,
and get smooth maps $\exp\colon {\cal D}_{\exp} \rightarrow M$ and $\exp_x:=\exp|_{T_xM\cap \calD_{\exp}} \colon {\cal D}_{\exp_x} \rightarrow M$. Each geodesic $\alpha\colon J \rightarrow TM$ with $\alpha^\prime(0)=v$ satisfies $\alpha(t)=\exp(tv)$.
A manifold with an affine connection is called \emph{geodesically complete} if the domain of each maximal geodesic is all of $\RR$.

Let $V\subseteq {\cal D}_{\exp_x}$ be an open neighborhood of $0$ in $T_xM=:E$ that is star-shaped with respect to $0$ (i.e., $[0,1]V\subseteq V$) such that $\exp_x$ induces a diffeomorphism of $V$ onto its open\linebreak image $W$. Then $W$ is said to be a \emph{normal neighborhood of $x$}. We call the chart\linebreak $\varphi:=(\exp|_V^W)^{-1}\colon W \rightarrow V\subseteq E$ a \emph{normal chart at $x$}.
Normal neighborhoods do exist, as $\exp_x\colon {\cal D}_{\exp_x}\rightarrow M$ induces a local diffeomorphism at $0\in T_xM$, since $T_0\exp_x=\id_{T_xM}$ (cf.\ \cite[Th.~IV.4.1]{Lan01}).
 
%
%
The definition of $\nabla$ can be extended to tensor fields
of type $\lambda\colon(E_1,E_2)\mapsto \calL^n(E_1,E_2)$,\footnote{This means a tensor field of contravariant degree $1$ and covariant degree $n$.}
where $\calL^n(E_1,E_2)$ denotes the space of $n$-linear continuous maps $E_1^n\rightarrow E_2$. There is a unique collection
$$(\nabla^U\colon{\cal V}(U) \times \Gamma \calL^n(TU,TU) \rightarrow \Gamma \calL^n(TU,TU))_{U \subseteq M \, \mbox{\scriptsize open}}$$
of maps such that
$((\nabla^{U_1})_\xi\omega)|_{U_2} = (\nabla^{U_2})_{\xi|_{U_2}}\omega|_{U_2}$
for all $\xi\in{\cal V}(U_1)$, $\omega\in \Gamma \calL^n(TU_1,TU_1)$, $U_2\subseteq U_1\subseteq M$, and such that
$$((\nabla^U)_\xi\omega)(\eta_1,\dots,\eta_n)
	\ =\ (\nabla^U)_\xi(\omega(\eta_1,\dots,\eta_n))
	- \sum_{i=1}^n \omega(\eta_1,\dots,(\nabla^U)_\xi\eta_i,\dots,\eta_n),$$
where $\eta_1,\dots,\eta_n$ denote any appropriate smooth vector fields. That is why we have a derivation property with respect to the $n+1$ variables $\omega,\eta_1,\dots,\eta_n$. There will be no confusion when we shall often suppress the index set $U$.
As far as the vector field $\xi$ is concerned, $(\nabla\!_\xi\omega)(x)$ depends only on $\xi(x)$. Therefore it makes sense to define 
$\nabla\!_v\omega$ for vectors $v\in T_xM$, $x\in M$. Then we have
$$(\nabla\!_v\omega)(\eta_1(x),\dots,\eta_n(x))
	\ =\ \nabla\!_v (\omega(\eta_1,\dots,\eta_n))
	- \sum_{i=1}^n \omega(x)(\eta_1(x),\dots,\nabla\!_v\eta_i,\dots,\eta_n(x)).$$

The definition of $\nabla\!_{\alpha^\prime}$ can be extended to lifts of $\alpha$ into $\calL^n(TM,TM)$: Given a curve $\alpha\colon J\rightarrow M$ and lifts $\omega\colon J\rightarrow \calL^n(TM,TM)$ and $\gamma_1,\dots,\gamma_n\colon J \rightarrow TM$ of $\alpha$, we denote by $\omega(\gamma_1,\dots,\gamma_n)$ the lift of $\alpha$ to $TM$ defined by
$$\omega(\gamma_1,\dots,\gamma_n)(t)
	\ =\ \omega(t)(\gamma_1(t),\dots,\gamma_n(t)).$$
The \emph{derivative of a lift $\omega$ along $\alpha$} is the unique lift $\nabla\!_{\alpha^\prime}\omega$ of $\alpha$ to $\calL^n(TM,TM)$ that satisfies
$$(\nabla\!_{\alpha^\prime}\omega)(\gamma_1,\dots,\gamma_n)
	\ =\ \nabla\!_{\alpha^\prime} (\omega(\gamma_1,\dots,\gamma_n))
	- \sum_{i=1}^n \omega(
		\gamma_1,\dots,\nabla\!_{\alpha^\prime}\gamma_i,\dots,\gamma_n)$$
for all lifts $\gamma_1,\dots,\gamma_n$ of $\alpha$ to $TM$.
We also use the notation $\nabla\!_{\alpha^\prime(t)}\omega$. Then we have
$$(\nabla\!_{\alpha^\prime(t)}\omega)(\gamma_1(t),\dots,\gamma_n(t))
	\ =\ \nabla\!_{\alpha^\prime(t)} (\omega(\gamma_1,\dots,\gamma_n))
	- \sum_{i=1}^n \omega(t)(\gamma_1(t),\dots,\nabla\!_{\alpha^\prime(t)}\gamma_i,\dots,\gamma_n(t)).$$
A lift $\omega$ of $\alpha$ is said to be {\em $\alpha$-parallel} if $\nabla\!_{\alpha^\prime}\omega = 0$. In this case, parallel transport along $\alpha$ commutes with $\omega$ in the sense that
$$\omega(t_1)\big(P_{t_0}^{t_1}(\alpha)(v_1),\ldots,P_{t_0}^{t_1}(\alpha)(v_n)\big) \ =\ P_{t_0}^{t_1}(\alpha)\big(\omega(t_0)(v_1,\ldots,v_n)\big)$$
for all $t_0,t_1\in J$ and $v_1,\ldots,v_n\in T_{\alpha(t_0)}M$.

The following proposition brings the definitions together.
\begin{proposition} \label{prop:derivativesAlongCurves}
	Let $\alpha\colon J \rightarrow M$ be a curve on a manifold $M$. Then we have:
	\begin{enumerate}
		\item[\rm (1)] For every smooth vector field $\eta$ on $M$, the lift $\eta\circ\alpha\colon J \rightarrow TM$ of $\alpha$ satisfies
			$$\nabla\!_{\alpha^\prime(t)}(\eta\circ\alpha)
				\ =\ \nabla\!_{\alpha^\prime(t)}\eta$$
			for all $t\in J$.
		\item[\rm (2)] For every tensor field $\omega\colon M \rightarrow \calL^n(TM,TM)$ on $M$, the lift $\omega\circ\alpha\colon J \rightarrow \calL^n(TM,TM)$ of $\alpha$ satisfies
			$$\nabla\!_{\alpha^\prime(t)}(\omega\circ\alpha)
				\ =\ \nabla\!_{\alpha^\prime(t)}\omega$$
			for all $t\in J$.
	\end{enumerate}
\end{proposition}
\begin{proof}
	Cf.\ \cite[Cor.~VIII.3.2 and Cor.~VIII.3.6]{Lan01}.
\end{proof}

Further details can be found in \cite[IV, VIII and X]{Lan01}, but basically for the case of torsionfree connections.
Cf.\ also \cite{KN63}, \cite{Kli82} and \cite{Ber08} for more material on connections.

\subsection{Affine Maps}
\label{sec:affineMaps}
Given two affine manifolds $(M_1,B_1)$ and $(M_2,B_2)$, a map $f\colon M_1\rightarrow M_2$ is called \emph{affine}, if $TTf\circ B_1=B_2\circ (Tf\oplus Tf)$. Working with charts $\varphi_1\colon U_1\rightarrow V_1\subseteq E_1$ of $M_1$ and\linebreak $\varphi_2\colon U_2\rightarrow V_2\subseteq E_2$ of $M_2$ such that $f(U_1)=U_2$, this can be written as
\begin{equation}\label{eqn:affineMaps}
	d^2f^\varphi(x)(v,w)+df^\varphi(x)((B_1^{\varphi_1})_x(v,w))\ =\ (B_2^{\varphi_2})_{f^\varphi(x)}(df^\varphi(x)(v),df^\varphi(x)(w))
\end{equation}
for all $x$ in the domain of the local representation $f^\varphi\colon V_1\rightarrow V_2$ of $f$ and $v,w\in E_1$.

Affine maps are compatible with parallel transport along curves, i.e., $T_{\alpha(t_1)}f\circ P_{t_0}^{t_1}(\alpha) = P_{t_0}^{t_1}(f\circ\alpha)\circ T_{\alpha(t_0)}f$ for all curves $\alpha\colon J\rightarrow M_1$ with $t_0,t_1\in J$. Geodesics  are mapped to geodesics. Further, we have $Tf(\calD_{\exp,1})\subseteq \calD_{\exp,2}$ and $f\circ\exp=\exp\circ Tf$.
A consequence is that, given an affine map, its values on connected components are uniquely determined by the tangent map at a single point, i.e., given affine maps $f,g\colon M_1\rightarrow M_2$ with $T_xf=T_xg$ for some $x\in M_1$, we have $f=g$ if $M_1$ is connected (cf.\ proof of \cite[Lem.~3.5]{Nee02}).

Affine maps are compatible with covariant derivatives of related vector fields, i.e., \linebreak $Tf(\nabla\!_v\eta_1)=\nabla\!_{Tf(v)}\eta_2$ for all $v\in TM_1$ and $\eta_1\in\calV(M_1)$, $\eta_2\in\calV(M_2)$ with $Tf\circ\eta_1=\eta_2\circ f$.

\subsection{Affine and Infinitesimal Affine Automorphisms}
\label{sec:affAndInfAffAuto}
Let $(M,\nabla)$ ($=(M,B)$) be an affine Banach manifold.
A diffeomorphism $f\in\Diff(M)$ is called an \emph{affine automorphism} if it is affine.
A vector field $\xi\in\calV(M)$ is called an \emph{infinitesimal affine automorphism} if each flow map $\flow^\xi_t$ is an affine isomorphism. We denote the set of all affine automorphisms by $\Aut(M,\nabla)$ and the set of all infinitesimal automorphisms by $\Kill(M,\nabla)$ or $\Kill(M,B)$.

The property of a vector field $\xi\in\calV(M)$ to be an infinitesimal affine automorphism can be checked locally: For a chart $\varphi\colon U\rightarrow V\subseteq E$ of $M$, the local representation $\overline\xi^\varphi\colon V\rightarrow E$ must satisfy
\begin{eqnarray}
	\lefteqn{d^2\xi^\varphi(x)(v,w)+d\xi^\varphi(B^\varphi_x(v,w))} \nonumber \\
	&=& dB^\varphi(x)(\xi^\varphi(x))(v,w)+B^\varphi_x(d\xi^\varphi(x)(v),w) + B^\varphi_x(v,d\xi^\varphi(x)(w)) \label{eqn:infAffAutoChart}
\end{eqnarray}
for all $x\in V$ and $v,w\in E$ (cf.\ \cite[Rem.~3.10]{Klo09a}).

We assume that $M$ is pure, i.e., that it has a single model space $E$.
Then the set
$\Fr(M):=\cup_{x\in M}\Iso(E,T_xM)$
(of topological linear isomorphisms) equipped with the projection $q\colon \Fr(M)\rightarrow M,\ \Iso(E,T_xM)\ni p \mapsto x$ carries the structure of a smooth $\GL(E)$-principal bundle with respect to the action
$$\rho\colon\Fr(M) \times \GL(E)\rightarrow \Fr(M),\ (p,g)\mapsto p.g:=p\circ g.$$

More precisely, for each chart $\varphi\colon U\rightarrow V\subseteq E$ of $M$, the map
$$\begin{array}{ccccc}
	\Fr(\varphi)\colon & \Fr(U) &\rightarrow & V\times \GL(E)& \subseteq E\times\gl(E) \\
	&\Iso(E,T_xU)\ni p & \mapsto & (\varphi(x),d\varphi(x)\circ p).
\end{array}$$
is a bundle chart of $\Fr(M)$, and we have $$q(\Fr(\varphi)^{-1}(\varphi(x),g))=x \quad\mbox{and}\quad \Fr(\varphi)^{-1}(\varphi(x),g_1g_2)=\Fr(\varphi)^{-1}(\varphi(x),g_1).g_2$$
for all $x\in U$ and $g, g_1,g_2\in \GL(E)$. The bundle $\Fr(M)$ is called the \emph{frame bundle over $M$}. For further details, see \cite[7.10.1]{Bou07}.

A diffeomorphism $f\colon M_1\rightarrow M_2$ of Banach manifolds (modelled on $E$) induces a principal bundle isomorphism $\Fr(f)\colon \Fr(M_1)\rightarrow \Fr(M_2)$ over $f$ defined by $\Fr(f)(p)=T_xf\circ p$ where $p\in \Iso(E,T_xM_1)$.

\begin{proposition}[{cf.\ \cite[Prop.~3.11 and Cor.~3.12]{Klo09a}}]\label{prop:isomorphismOfKill}
	Given an affine Banach manifold $(M,\nabla)$, the set $\Kill(M,\nabla)$ is a subalgebra of the Lie algebra $\calV(M)$ of smooth vector fields on $M$. If $M$ is connected, $\Kill(M,\nabla)$ can be turned into a Banach--Lie algebra whose Banach space structure is uniquely determined by the requirement that for each $p\in\Fr(M)$, the map
	$$\Kill(M,\nabla)\rightarrow T_p(\Fr(M)),\ \xi\mapsto \left.\frac{d}{dt}\right|_{t=0}\Fr(\flow^{\xi}_t)(p)$$
	is a closed embedding. Then, for each $x\in M$, also the map
	$$\Kill(M,\nabla)\rightarrow T_xM \times \calL(T_xM,T_xM),\ \xi\mapsto \big(\xi(x),\ v\mapsto\nabla\!_v\xi\big)$$
	is a closed embedding of Banach spaces.
\end{proposition}
\begin{theorem}[{cf.\ \cite[Th.~3.14]{Klo09a}}]\label{th:Kill(M,nabla)Complete}
	Let $(M,\nabla)$ be an affine Banach manifold. If it is geodesically complete, then all vector fields in $\Kill(M,\nabla)$ are complete.
\end{theorem}
\begin{theorem}[{cf.\ \cite[Th.~3.15]{Klo09a}}]\label{th:autoGroupOfAffineManifold}
	Let $(M,\nabla)$ be a connected affine Banach manifold that is geodesically complete.
	The automorphism group $\Aut(M,\nabla)$ can be turned into a Banach--Lie group such that
	$$\exp\colon \Kill(M,\nabla)\rightarrow \Aut(M,\nabla),\ \xi \mapsto \flow^{-\xi}_1$$
	is its exponential map. The natural map $\tau\colon \Aut(M,\nabla)\times M \rightarrow M$ is a smooth action whose derived action is the inclusion map $\Kill(M,\nabla)\hookrightarrow\calV(M)$, i.e., $-T\tau(\id_M,x)(\xi,0)=\xi(x)$.
\end{theorem}
%
%
%
%
%
%
%
\section{Local Integration of Maps between Tangent Spaces to Affine Maps}\label{sec:locIntegration}
Let $(M_1,\nabla\!_1)$ and $(M_2,\nabla\!_2)$ be affine Banach manifolds with distinguished points $b_1\in M_1$ and $b_2\in M_2$ called the \emph{base points}. A continuous linear map $A\colon T_{b_1}M_1\rightarrow T_{b_2}M_2$ between the tangent spaces at the base points is called \emph{locally integrable} if there exists an affine map $f$ from an open neighborhood of $b_1$ into $M_2$ that satisfies $T_{b_1}f=A$.

The main result of this section is that for affine manifolds with parallel torsion and curvature tensors, such a map $A$ is locally integrable if it intertwines these tensors in the base points. The proof given here is close to that of O.~Loos for the finite-dimensional case (cf.\ \cite[pp.~104-111]{Loo69}).

\subsection{Torsion and Curvature Tensor}
\label{sec:torsionAndCurvature}
Let $(M,\nabla)$ be an affine Banach manifold.
There exists a unique tensor field $\Tor$ of type $\lambda:(E_1,E_2)\mapsto \calL^2(E_1,E_2)$ on $M$ such that for any open set $U$ in $M$ and smooth vector fields $\xi,\eta$ on $U$, we have
$$\Tor(\xi,\eta) \ =\ \nabla\!_\xi\eta-\nabla\!_\eta\xi-[\xi,\eta].$$
We shall call $\Tor$ the {\em torsion tensor}. Given a chart $\varphi\colon U \rightarrow V \subseteq E$, a point $x\in V$ and vectors $v,w\in E$, we have
$$\Tor^\varphi_x(v,w) \ =\ B^\varphi_x(v,w) -B^\varphi_x(w,v).$$
The torsion tensor $\Tor$ is skew-symmetric, i.e.,
$\Tor_x(v,w)=-\Tor_x(w,v)$
for all $x\in M$ and $v,w\in T_xM$. If it is identically zero, then the affine connection is called \emph{torsionfree} or also \emph{symmetric} (and can then be described by using a spray, cf.\ \cite[IV, \S 3]{Lan01} or \cite[III.11]{Ber08}).

There exists a unique tensor field $R$ of type $\lambda:(E_1,E_2)\mapsto \calL^3(E_1,E_2)$ on $M$ such that for any open set $U$ in $M$ and smooth vector fields $\xi,\eta,\zeta$ on $U$, we have
$$R(\xi,\eta,\zeta) \ =\ \nabla\!_\xi\nabla\!_\eta\zeta - \nabla\!_\eta\nabla\!_\xi\zeta - \nabla\!_{[\xi,\eta]}\zeta.$$ 
We shall call $R$ the {\em curvature tensor}.
Given a chart $\varphi\colon U \rightarrow V \subseteq E$, a point $x\in V$ and vectors $v,w,z\in E$, we have
\begin{eqnarray*}
	R^\varphi_x(v,w,z) 
		& = & B^\varphi_x(B^\varphi_x(z,w),v) - B^\varphi_x(B^\varphi_x(z,v),w)\\
	& & {}+ dB^\varphi(x)(w)(z,v) - dB^\varphi(x)(v)(z,w).
\end{eqnarray*}
The curvature tensor $R$ is skew-symmetric with respect to the first two arguments, i.e.,
$R_x(v,w,z)=-R_x(w,v,z)$
for all $x\in M$ and $v,w,z\in T_xM$.
For further details, see \cite[IX, \S 1]{Lan01}.

Every affine map $f$ between two affine Banach manifolds $(M_1,\nabla\!_1)$ and $(M_2,\nabla\!_2)$, intertwines the torsion and curvature tensors in the sense that
$$T_xf\circ (\Tor_1)_x = (\Tor_2)_{f(x)}\circ (T_xf)^2 \quad\mbox{and}\quad T_xf \circ (R_1)_x = (R_2)_{f(x)}\circ (T_xf)^3$$
for all $x\in M_1$. This can be checked by working with charts.
\subsection{Exterior Derivatives and Wedge Products}
\label{sec:extDerivAndWedgeProd}
Let $M$ be a Banach manifold and $F$ a Banach space.
The {\em exterior derivative} of an $F$-valued $n$-form $\omega$ on $M$ is the unique $F$-valued $(n+1)$-form $\extd\omega$ on $M$ such that for any open set $U$ in $M$ and smooth vector fields $\xi_0,\dots,\xi_n$ on $U$, we have
\begin{eqnarray*}
	\lefteqn{\extd\omega(\xi_0, \dots, \xi_n)} \\
	& \quad = & \sum_{i=0}^n (-1)^i \xi_i.(\omega
			(\xi_0,\dots,\widehat\xi_i,\dots,\xi_n)) \\
	& & {} + \sum_{i<j} (-1)^{i+j} \omega([\xi_i,\xi_j],\xi_0,\dots,
				\widehat\xi_i, \dots, \widehat\xi_j, \dots, \xi_n).
\end{eqnarray*}
Given a chart $\varphi\colon U \rightarrow V \subseteq E$, its local representation $(\extd\omega)^\varphi$ is given by
$$(\extd\omega)^\varphi_x(v_0, \dots, v_n)
	= \sum_{i=0}^n (-1)^i d\omega^\varphi(x)(v_i)(v_0, \dots,\widehat{v_i}, \dots, v_n).$$
with $x\in V$ and vectors $v_0,\dots,v_n \in E$.\footnote{The local representation $\omega^\varphi\colon V\rightarrow \Alt^n(E,F)$ of the $n$-form $\omega$ can be considered as a Banach space valued map with derivative $d\omega^\varphi\colon V\rightarrow \calL(E,\Alt^n(E,F))$.}
For every form $\omega$, we have $\extd\extd\omega = 0$.

Let $\omega_1$ and $\omega_2$ be differential forms of degree $n_1$ and $n_2$ on $M$ with values in $F_1$ and $F_2$, respectively, and let $\beta\colon F_1\times F_2\rightarrow F$ be a continuous bilinear map into a Banach space $F$. The {\em wedge product $\omega_1 \wedge_\beta \omega_2$ of $\omega_1$ and $\omega_2$} is defined as the $F$-valued $(n_1+n_2)$-form on $M$ given by
\begin{eqnarray*}
	\lefteqn{(\omega_1 \wedge_\beta \omega_2)_x(v_1,\dots,v_{n_1+n_2})} \\
	& \quad = & \frac{1}{n_1!\,n_2!}\sum_{\sigma\in S_{n_1+n_2}} \sgn(\sigma)
			\beta\big((\omega_1)_x(v_{\sigma(1)},\dots,v_{\sigma(n_1)}),
			(\omega_2)_x(v_{\sigma(n_1+1)},\dots,v_{\sigma(n_1+n_2)})\big).
\end{eqnarray*}
For its exterior derivative, we have the formula
$$\extd(\omega_1 \wedge_\beta \omega_2)
	= \extd\omega_1 \wedge_\beta \omega_2 + (-1)^{n_1}\omega_1\wedge_\beta \extd\omega_2.$$

Pull-backs are compatible with exterior derivatives and wedge products, i.e., given a smooth map $f\colon N\rightarrow M$ of manifolds, we have formulas like
$$\extd(f^\ast\omega)=f^\ast(\extd\omega) \quad \mbox{and} \quad
	f^\ast(\omega_1\wedge_\beta\omega_2)
	= f^\ast\omega_1\wedge_\beta f^\ast\omega_2.$$
For further details, cf.\ \cite[V, \S 3]{Lan01} concerning real-valued differential forms where $\beta$ is the multiplication in $\RR$, and cf.\ \cite[8.3-8.5]{Bou07} concerning the general case.
\subsection{Structure Equations}
\label{sec:structureEquations}
Throughout this section, let $(M,\nabla,b)$ be an affine Banach manifold $M$ with a point $b\in M$ called the {\em base point}. Let $\varphi\colon W \rightarrow V\subseteq T_bM =: E$ be a normal chart at $b$.
The base point $b$ can be joined with each $x\in W$ by a unique\footnote{Indeed, given a geodesic
	$\alpha\colon [0,1]\rightarrow W$ with $\alpha(t):=\exp_b(tv)$ and $v\in T_bM$ that joins $b$ with $x=\exp_b(v)$, its local representation $\alpha^\varphi$ maps the interval $[0,1]$ onto a compact set in $V$, so that there is a $\lambda\in {[0,1[}$ such that $\alpha^\varphi([0,1])\subseteq \lambda V$. If $v \notin V$, then there is a $t\in[0,1]$ with $tv \in V\backslash \lambda V$, so that $\alpha^\varphi(t)=\varphi(\exp_b(tv)) = tv \in V\backslash\lambda V$ contradicts $\alpha^\varphi([0,1])\subseteq \lambda V$. Hence $v$ lies in $V$, so that $v=\varphi(x)$.}
geodesic $\alpha_x\colon [0,1]\rightarrow W$. It is given by $\alpha_x(t):=\exp_b(t\varphi(x))$.
Consequently, the exponential $\exp_b^W$ of the open submanifold $W$ is given by $\exp_b|_V$.

Given a vector $v\in T_bM$, we define the smooth vector field
$$v^\ast\colon W \rightarrow TM,\ x \mapsto P_0^1(\alpha_x)(v)$$
and call it an {\em adapted vector field}. Note that, for each $x\in M$, the tangent space $T_xM$ is given by $\{v^\ast(x)\colon v\in T_bM\}$.
For working in the chart $\varphi$, we define
$\bar{v}^\ast\colon V \rightarrow  E$
for each $\bar{v}\in E$ such that $\bar{v}^\ast$ is the local representation of $v^\ast$ where $v:=T\varphi^{-1}(0,\bar{v})$.
We shall define several differential forms on $W$ that describe the affine structure.
The correctness of the following definitions can be checked by working with local formulas.
The map $\theta$ is defined as the $T_bM$-valued $1$-form on $W$ given by
$$\theta_x=P_1^0(\alpha_x) \colon T_xW\rightarrow T_bM,$$
and the \emph{connection form} $\omega$ is defined as the $\calL(T_bM,T_bM)$-valued $1$-form on $W$ given by
$$\omega_x\colon T_xW\rightarrow \calL(T_bM,T_bM),\ \omega_x(v)(w):=\theta_x(\nabla\!_v w^\ast)$$
with local representation
$$\omega^\varphi_{\bar{x}}\colon E \rightarrow \calL(E,E),\ \omega^\varphi_{\bar{x}}(\bar{v})(\bar{w}) := \theta^\varphi_{\bar{x}}\big(d\bar{w}^\ast(\bar{x})(\bar{v}) - B^\varphi_{\bar{x}}(\bar{w}^\ast(\bar{x}),\bar{v})\big).$$
The \emph{torsion form} $\Theta$ is defined as the $T_bM$-valued 2-form on $W$ given by
$$\Theta_x\colon T_xW\times T_xW \rightarrow T_bM,\ \Theta_x(v,w)=\theta_x(\Tor_x(v,w)),$$
and the \emph{curvature form} $\Omega$ is defined as the $\calL(T_bM,T_bM)$-valued $2$-form given by
$$\Omega_x\colon T_xW\times T_xW \rightarrow \calL(T_bM,T_bM),\  \Omega_x(v,w)(z) := \theta_x\big(R_x(v,w,z^\ast(x))\big).$$
Notice that the alternating properties are due to $\Tor_x(v,v)=0$ and $R_x(v,v,\cdot)=0$, respectively.
Now, we can formulate the structure equations that describe the affine connection.
\begin{proposition}[Structure equations of \'E.~Cartan] \label{prop:structureEquations}
	Let $\ev$ denote the evaluation map
	$$\ev\colon \calL(T_bM,T_bM)\times T_bM \rightarrow T_bM,\ (A,v) \mapsto A(v)$$
	and $\Gamma$ the composition map
	$$\Gamma\colon \calL(T_bM,T_bM)\times \calL(T_bM,T_bM) \rightarrow \calL(T_bM,T_bM),\ (A,B) \mapsto A\circ B.$$
	Then the following equations hold:
	\begin{enumerate}
		\item[\rm (1)] $\extd\theta + \omega\wedge_{\ev}\theta = \Theta$.
		\item[\rm (2)] $\extd\omega + \omega\wedge_\Gamma\omega = \Omega.$
	\end{enumerate}
\end{proposition}
\begin{proof}
	For all $v\in T_bM$, we put $\ev_v:=\ev(\cdot,v)$. In the following, given any maps\linebreak $f,g\colon W\rightarrow \calL(T_bM,T_bM)$ and $h\colon W\rightarrow T_bM$, we denote by $\ev(f,h)\colon W\rightarrow T_bM$,\linebreak $\ev_v(f)\colon W\rightarrow T_bM$ and $\Gamma(f,g)\colon W\rightarrow \calL(T_bM,T_bM)$ the maps defined by $\ev(f,h)(x):=\ev(f(x),h(x))$, $\ev_v(f)(x):=\ev_v(f(x))$  and $\Gamma(f,g)(x):=\Gamma(f(x),g(x))$, respectively.
	
	(1) The proof works as in the finite-dimensional case (cf.\ \cite[p.~106]{Loo69}).
	It suffices to show that, given any vectors $v,w \in T_bM$, we have
	$$\extd\theta(v^\ast,w^\ast) + (\omega\wedge\theta)(v^\ast,w^\ast) \ =\ \Theta(v^\ast,w^\ast).$$
	We observe
	\begin{eqnarray*}
		\lefteqn{\extd\theta(v^\ast,w^\ast) + (\omega\wedge\theta)(v^\ast,w^\ast)} \\
		&\quad =& (v^\ast).(\theta(w^\ast)) - (w^\ast).(\theta(v^\ast)) - \theta([v^\ast,w^\ast]) + 
		\ev(\omega(v^\ast),\theta(w^\ast)) - \ev(\omega(w^\ast),\theta(v^\ast)) \\
		&\quad =& 0 - 0 - \theta([v^\ast,w^\ast]) + \theta(\nabla\!_{v^\ast}w^\ast) - \theta(\nabla\!_{w^\ast}v^\ast)
		\ = \ \theta (-[v^\ast,w^\ast] + \nabla\!_{v^\ast}w^\ast - \nabla\!_{w^\ast}v^\ast) \\
		&\quad =& \theta(\Tor(v^\ast,w^\ast))
		\ = \ \Theta(v^\ast,w^\ast),
	\end{eqnarray*}
	since $\theta(w^\ast)(x)\equiv w$ and $\theta(v^\ast)(x)\equiv v$.
	
	(2) Cf.\ \cite[pp.~106, 107]{Loo69} for the finite-dimensional case.
	There, vector fields on $W$ can be represented by $\calF(W)$-linear combinations of adapted vector fields, where $\calF(W)$ denotes the algebra of real functions on $W$.	
	
	It suffices to show that, given any vectors $v,w,z \in T_bM$, we have
	$$\ev_z(\extd\omega(v^\ast,w^\ast)) + \ev_z((\omega\wedge\omega)(v^\ast,w^\ast))
		\ =\ \ev_z(\Omega(v^\ast,w^\ast)).$$
	We observe
	\begin{eqnarray*}
		\ev_z(\extd\omega(v^\ast,w^\ast))
		&=& \ev_z\big((v^\ast).(\omega(w^\ast))\big) - \ev_z\big((w^\ast).(\omega(v^\ast))\big) - \ev_z(\omega([v^\ast,w^\ast])) \\
		&=& \ev_z\big((v^\ast).(\omega(w^\ast))\big) - \ev_z\big((w^\ast).(\omega(v^\ast))\big) - \theta(\nabla\!_{[v^\ast,w^\ast]}z^\ast)
	\end{eqnarray*}
	and
	$$\ev_z(\Omega(v^\ast,w^\ast))
		\ =\ \theta(R(v^\ast,w^\ast,z^\ast))
		\ =\ \theta(\nabla\!_{v^\ast}\nabla\!_{w^\ast}z^\ast - \nabla\!_{w^\ast}\nabla\!_{v^\ast}z^\ast - \nabla\!_{[v^\ast,w^\ast]}z^\ast),$$
	so that we have to show
	\begin{eqnarray*}
		&& \ev_z\big((v^\ast).(\omega(w^\ast))\big) - \ev_z\big((w^\ast).(\omega(v^\ast))\big) + \ev_z\big((\omega\wedge\omega)(v^\ast,w^\ast)\big) \\
		&& \qquad = \  \theta(\nabla\!_{v^\ast}\nabla\!_{w^\ast}z^\ast - \nabla\!_{w^\ast}\nabla\!_{v^\ast}z^\ast).
	\end{eqnarray*}
	Since
	$$\ev_z((\omega\wedge\omega)(v^\ast,w^\ast)) \ =\ \ev_z\big(\Gamma(\omega(v^\ast),\omega(w^\ast))-\Gamma(\omega(w^\ast),\omega(v^\ast))\big),$$
	it suffices to show
	$$\ev_z\big((v^\ast).(\omega(w^\ast))\big) + \ev_z\big(\Gamma(\omega(v^\ast),\omega(w^\ast))\big) \ =\ \theta(\nabla\!_{v^\ast}\nabla\!_{w^\ast}z^\ast),$$
	because then the same equation holds if $v$ and $w$ are interchanged.
	We have
	$$\ev_z\big((v^\ast).(\omega(w^\ast))\big)
		\ =\ \ev_z(d(\omega(w^\ast))\circ v^\ast)
		\ =\ d\big(\ev_z(\omega(w^\ast))\big)\circ v^\ast
		\ =\ d(\theta(\nabla\!_{w^\ast}z^\ast)) \circ v^\ast$$
	and
	$$\ev_z\big(\Gamma(\omega(v^\ast),\omega(w^\ast))\big) \ =\ \ev(\omega(v^\ast),\theta(\nabla\!_{w^\ast}z^\ast)).$$
	Abbreviating the vector field $\nabla\!_{w^\ast}z^\ast$ by $\xi$, we have to prove
	$$d(\theta(\xi)) \circ v^\ast + \ev(\omega(v^\ast),\theta(\xi)) \ =\ \theta(\nabla\!_{v^\ast}\xi).$$
	Working in the chart $\varphi$, we shall show
	$$d(\theta^\varphi(\xi^\varphi))(\bar x)(v_{\bar x})
		+ \omega^\varphi_{\bar x}(v_{\bar x})\big(\theta^\varphi_{\bar x}(\xi^\varphi(\bar x))\big)
		\ =\ \theta^\varphi_{\bar x}\big(d\xi^\varphi(v_{\bar x})-B^\varphi_{\bar x}(\xi^\varphi(\bar x),v_{\bar x})\big)$$
	for all $\bar x\in V$, where $v_{\bar x}:=(v^\ast)^\varphi(\bar x)$.
	Putting $\bar u:= \theta^\varphi_{\bar x}(\xi^\varphi(\bar x))$, we have $\xi^\varphi(\bar x) = \bar u^\ast(\bar x)$.
	Computing the left hand side by applying the product rule to the first summand and the formula for $\omega^\varphi_{\bar x}$ to the second one, we obtain
	$$(d\theta^\varphi(\bar x)(v_{\bar x}))(\xi^\varphi(\bar x)) + \theta^\varphi_{\bar x}(d\xi^\varphi(\bar x)(v_{\bar x}))
	+ \theta^\varphi_{\bar x}\big(du^\ast(\bar{x})(v_{\bar x}) - B^\varphi_{\bar x}(u^\ast(\bar x),v_{\bar x})\big).$$
	Therefore, we merely have to verify
	$$(d\theta^\varphi(\bar x)(v_{\bar x}))(\bar u^\ast(\bar x))
		+ \theta^\varphi_{\bar x}(du^\ast(\bar{x})(v_{\bar x}))
		\ =\ 0.$$
	Applying the product rule backwards yields the equivalent equation
	$d(\theta^\varphi(\bar u^\ast))(\bar x)(v_{\bar x}) = 0,$
	which is true, since the differential vanishes, the map $\theta^\varphi(u^\ast)$ being constant.
\end{proof}
Now we put $\widehat V:=\{(t,v)\in\RR\times T_bM:tv \in V\}$ and consider the map
$\Phi\colon \widehat V \rightarrow W$ defined by $\Phi(t,v):=\exp_b(tv)$.
Its derivative is given by
\begin{equation} \label{eqn:TPhi}
	T_{(t,v)}\Phi(t^\prime,v^\prime)
	\ =\ T_{tv}\exp_b(t^\prime v + tv^\prime)
	\ =\ t^\prime T_{tv}\exp_b(v) + tT_{tv}\exp_b(v^\prime). 
\end{equation}
Given $(t,v)\in \widehat V$, we have
\begin{equation} \label{eqn:T_tvExp(v)1}
	v^\ast(\exp_b(tv)) \ =\ \frac{d}{dt}\exp_b(tv) \ =\ T_{tv}\exp_b(v)
\end{equation}
and hence
\begin{equation} \label{eqn:T_tvExp(v)2}
	\theta_{\exp_b(tv)}(T_{tv}\exp_b(v)) \ =\ v.
\end{equation}
Let us denote the projections from $\widehat V \subseteq \RR\times T_bM$ onto its components by $\lambda_{\RR}$ and $\lambda_{T_bM}$. Sometimes, we shall use these symbols also for the projections from all of $\RR\times T_bM$, but there will be no confusion.
\begin{lemma}
\label{lem:phi*thetaOmega}
	Let $\widehat\theta\colon \widehat V \rightarrow \calL(\RR\times T_bM, T_bM)$ be the $T_bM$-valued $1$-form defined by
	$$\widehat\theta_{(t,v)}(t^\prime,v^\prime)
		:=t\theta_{\exp_b(tv)}(T_{tv}\exp_b(v^\prime))$$
	and $\widehat\omega\colon \widehat V \rightarrow \calL(\RR\times T_bM,\calL(T_bM,T_bM))$ the $\calL(T_bM,T_bM)$-valued $1$-form defined by
	$$\widehat\omega_{(t,v)}(t^\prime,v^\prime)
		:=t\omega_{\exp_b(tv)}(T_{tv}\exp_b(v^\prime)).$$
	Then we have
	\begin{enumerate}
		\item[\rm(1)] $\Phi^\ast\theta = \lambda_{T_bM}\,d\lambda_{\RR} + \widehat\theta\colon \quad \widehat V \rightarrow \calL(\RR\times T_bM, T_bM)$ \quad and
		\item[\rm(2)] $\Phi^\ast\omega = \widehat\omega\colon \quad \widehat V \rightarrow \calL(\RR\times T_bM,\calL(T_bM,T_bM)).$
	\end{enumerate}
\end{lemma}
\begin{proof}
	The proof works as in the finite-dimensional case (cf.\ \cite[p.~107]{Loo69}).
	
	(1) By using (\ref{eqn:TPhi}) and (\ref{eqn:T_tvExp(v)2}), we get:
	\begin{eqnarray*}
		(\Phi^\ast\theta)_{(t,v)}(t^\prime,v^\prime)
			&=& \theta_{\Phi(t,v)}(T_{(t,v)}\Phi(t^\prime,v^\prime)) \\
		&=& \theta_{\exp_b(tv)}(t^\prime T_{tv}\exp_b(v) + tT_{tv}\exp_b(v^\prime)) \\
		&=& t^\prime\theta_{\exp_b(tv)}(T_{tv}\exp_b(v))
			+ t\theta_{\exp_b(tv)}(T_{tv}\exp_b(v^\prime)) \\
		&=& t^\prime v + \widehat\theta_{(t,v)}(t^\prime,v^\prime).
	\end{eqnarray*}
	
	(2) By a simple computation, as in the proof of (1), we get
	\begin{eqnarray*}
		(\Phi^\ast\omega)_{(t,v)}(t^\prime,v^\prime) & = & 
			\omega_{\Phi(t,v)}(T_{(t,v)}\Phi(t^\prime,v^\prime)) \\
		& = & t^\prime\omega_{\exp_b(tv)}(T_{tv}\exp_b(v))
			+ t\omega_{\exp_b(tv)}(T_{tv}\exp(v^\prime)) \\
		& = & t^\prime\omega_{\exp_b(tv)}(T_{tv}\exp_b(v))
			+ \widehat\omega_{(t,v)}(t^\prime,v^\prime).
	\end{eqnarray*}
	It remains to check that the first summand vanishes.
	For every $z\in T_bM$, we have
	\begin{eqnarray*}
	t^\prime\omega_{\exp_b(tv)}(T_{tv}\exp_b(v))(z) & = &
		t^\prime\theta_{\exp_b(tv)}(\nabla\!_{T_{tv}\exp_b(v)}z^\ast)\\
	& \stackrel{\mbox{\scriptsize (\ref{eqn:T_tvExp(v)1})}}{=} & t^\prime\theta_{\exp_b(tv)}(\nabla\!_{\alpha^\prime(t)}z^\ast)
	\end{eqnarray*}
	with the geodesic $\alpha\colon [0,t]\rightarrow W$ defined by $\alpha(t):=\exp_b(tv)$.
	It suffices to show	$\nabla\!_{\alpha^\prime(t)}z^\ast = 0$.
	This is true, as $z^\ast\circ\alpha$ is $\alpha$-parallel.
\end{proof}
\begin{remark} \label{rem:widehatThetaOmega}
	Evidently, we have
	$$\widehat\theta_{(t,v)}(t^\prime,0)= 0 \quad \mbox{and} \quad 
		\widehat\omega_{(t,v)}(t^\prime,0)= 0.$$
\end{remark}
\begin{lemma} \label{lem:partialWidehatThetaOmega}
	The partial derivatives $\partial_1\widehat\theta$ and $\partial_1\widehat\omega$ can directly be expressed in terms of the exterior derivatives $\extd\widehat\theta$ and $\extd\widehat\omega$, respectively. We have
	\begin{enumerate}
		\item[\rm (1)] $\partial_1\widehat\theta(t,v)(t^\prime,v^\prime)
			\ =\ (\extd\widehat\theta)_{(t,v)}((1,0),(t^\prime,v^\prime))$ \quad and
		\item[\rm (2)] $\partial_1\widehat\omega(t,v)(t^\prime,v^\prime)
			\ =\ (\extd\widehat\omega)_{(t,v)}((1,0),(t^\prime,v^\prime))$
	\end{enumerate}
	for all $(t,v)\in\widehat V$ and $(t^\prime,v^\prime)\in \RR\times T_bM$.
\end{lemma}
\begin{proof}
	We have
	$$(\extd\widehat\theta)_{(t,v)}((1,0),(t^\prime,v^\prime))
		\ =\ d\widehat\theta(t,v)(1,0)(t^\prime,v^\prime)
		- d\widehat\theta(t,v)(t^\prime,v^\prime)(1,0),\footnote{Note that the 1-form $\widehat\theta\colon\widehat V\rightarrow\calL(\RR\times T_bM,T_bM)$ can also be considered as a Banach space valued map with derivative $d\widehat\theta$. }$$
	Since
	$d\widehat\theta(t,v)(1,0) \ =\ \partial_1\widehat\theta(t,v),$
	the minuend equals
	$\partial_1\widehat\theta(t,v)(t^\prime,v^\prime)$. The subtrahend vanishes, as
	$$d\widehat\theta(t,v)(t^\prime,v^\prime)(1,0)
		= d(\ev_{(1,0)}\circ\widehat\theta)(t,v)(t^\prime,v^\prime)$$
	and $\ev_{(1,0)}\circ\widehat\theta=0$ by Remark~\ref{rem:widehatThetaOmega}, where $\ev_{(1,0)}\colon \calL(\RR\times T_bM, T_bM)  \rightarrow T_bM$ denotes the evaluation map $A\mapsto A(1,0)$. This proves (1). An analogous argument shows (2).
\end{proof}
\begin{proposition}
\label{prop:diffEqnWidehatThetaOmega}
	The forms $\widehat\theta$ and $\widehat\omega$ satisfy the system of ordinary differential equations
	\begin{eqnarray}
		\partial_1\widehat\theta(t,v) &=& \lambda_{T_bM} + \widehat\omega_{(t,v)}\cdot v + (\Phi^\ast\Theta)_{(t,v)}((1,0),\cdot)
			\label{eqn:diffEqnWidehatTheta} \\
		\partial_1\widehat\omega(t,v) &=& (\Phi^\ast\Omega)_{(t,v)}((1,0),\cdot)
			\label{eqn:diffEqnWidehatOmega}
	\end{eqnarray}
	with initial conditions $\widehat\theta(0,v)=0$ and $\widehat\omega(0,v)=0$, where $\widehat\omega_{(t,v)}\cdot v \in \calL(\RR \times T_bM,T_bM)$ denotes the map defined by
	$$(\widehat\omega_{(t,v)}\cdot v)(t^\prime,v^\prime)
		:=\widehat\omega_{(t,v)}(t^\prime,v^\prime)(v).$$
\end{proposition}
\begin{proof}
	Cf.\ \cite[p.~108]{Loo69} for the finite-dimensional case, where basic representations of differential forms are used.
	The main idea is to take the pull-back of the structure equations by $\Phi$. The initial conditions are obvious by the definition of $\widehat\theta$ and $\widehat\omega$ (cf.\ Lemma~\ref{lem:phi*thetaOmega}).
	
	(\ref{eqn:diffEqnWidehatTheta}): From Proposition~\ref{prop:structureEquations}(1) and Lemma~\ref{lem:phi*thetaOmega}(1), we obtain
	$$\extd(\lambda_{T_bM}\,d\lambda_{\RR} + \widehat\theta) + \widehat\omega \wedge_{\ev}
		(\lambda_{T_bM}\,d\lambda_{\RR} + \widehat\theta) \ =\ \Phi^\ast\Theta.$$
	Expressing $\lambda_{T_bM}\,d\lambda_{\RR}$ by the wedge product $\lambda_{T_bM} \wedge_m d\lambda_{\RR}$, where $m\colon T_bM \times \RR \rightarrow T_bM$ denotes the scalar multiplication, we have
	$$\extd(\lambda_{T_bM}\,d\lambda_{\RR}) \ =\
		d\lambda_{T_bM}\wedge_m d\lambda_{\RR}
		+ (-1)^0\lambda_{T_bM}\wedge_m \extd(d\lambda_{\RR}) \ =\ d\lambda_{T_bM}\wedge_m d\lambda_{\RR},\footnote{Note that the maps $\lambda_{\RR}$ and $\lambda_{T_bM}$ can be considered as Banach space valued maps as well as $0$-forms, but that the derivatives $d\lambda_{\RR}$ and $d\lambda_{T_bM}$ are equal to the exterior derivatives $\extd\lambda_{\RR}$ and $\extd\lambda_{T_bM}$. In particular, we have $\extd(d\lambda_{\RR})=\extd\extd\lambda_{\RR}=0.$}$$
	so that
	$$d\lambda_{T_bM}\wedge_m d\lambda_{\RR}
		+ \extd\widehat\theta + \widehat\omega \wedge_{\ev} (\lambda_{T_bM}\,d\lambda_{\RR})
		+ \widehat\omega \wedge_{\ev}\widehat\theta \ =\ \Phi^\ast\Theta.$$
	We shall evaluate this at $(t,v)((1,0),(t^\prime,v^\prime))$:
	The first summand yields
	\begin{eqnarray*}
		\lefteqn{(d\lambda_{T_bM}\wedge_m
			d\lambda_{\RR})(t,v)((1,0),(t^\prime,v^\prime))} \\
		&\quad =& d\lambda_{T_bM}(t,v)(1,0)d\lambda_{\RR}(t,v)(t^\prime,v^\prime)
			- d\lambda_{T_bM}(t,v)(t^\prime,v^\prime)d\lambda_{\RR}(t,v)(1,0) \\
		&\quad =& 0t^\prime - v^\prime1 
		\ =\ -v^\prime
	\end{eqnarray*}
	and the second one yields
	$\partial_1\widehat\theta(t,v)(t^\prime,v^\prime)$
	by Lemma~\ref{lem:partialWidehatThetaOmega}. The third summand yields
	\begin{eqnarray*}
		\lefteqn{(\widehat\omega \wedge_{\ev}
			(\lambda_{T_bM}\,d\lambda_{\RR}))(t,v)((1,0),(t^\prime,v^\prime))} \\
		&\quad =& \widehat\omega_{(t,v)}(1,0)
			((\lambda_{T_bM}\,d\lambda_{\RR})(t,v)(t^\prime,v^\prime))
			- \widehat\omega_{(t,v)}(t^\prime,v^\prime)
			((\lambda_{T_bM}\,d\lambda_{\RR})(t,v)(1,0)) \\
		&\quad =& - \widehat\omega_{(t,v)}(t^\prime,v^\prime)(v),
	\end{eqnarray*}
	since $\widehat\omega_{(t,v)}(1,0)=0$ by Remark~\ref{rem:widehatThetaOmega},
	and the last one yields
	$$(\widehat\omega \wedge_{\ev}\widehat\theta)
			(t,v)((1,0),(t^\prime,v^\prime))
		\ = \ \widehat\omega_{(t,v)}(1,0)(\widehat\theta_{(t,v)}(t^\prime,v^\prime))
			- \widehat\omega_{(t,v)}(t^\prime,v^\prime)(\widehat\theta_{(t,v)}(1,0)) \ =\ 0, $$
	since $\widehat\omega_{(t,v)}(1,0)=0$ and $\widehat\theta_{(t,v)}(1,0)=0$ by Remark~\ref{rem:widehatThetaOmega}.
	We thus arrive at
	$$-v^\prime + \partial_1\widehat\theta(t,v)(t^\prime,v^\prime)
		- \widehat\omega_{(t,v)}(t^\prime,v^\prime)(v) \ =\ (\Phi^\ast\Theta)_{(t,v)}((1,0),(t^\prime,v^\prime)),$$
	which entails (\ref{eqn:diffEqnWidehatTheta}).
	
	(\ref{eqn:diffEqnWidehatOmega}): From Proposition~\ref{prop:structureEquations}(1) and Lemma~\ref{lem:phi*thetaOmega}(1), we obtain
	$$d\widehat\omega + \widehat\omega\wedge_\Gamma\widehat\omega
		\ =\ \Phi^\ast\Omega.$$
	We shall evaluate this at $(t,v)((1,0),(t^\prime,v^\prime))$:
	The first summand yields
	$\partial_1\widehat\omega(t,v)(t^\prime,v^\prime)$
	by Lemma~\ref{lem:partialWidehatThetaOmega}.
	The second one yields
	$$(\widehat\omega\wedge_\Gamma\widehat\omega)
			(t,v)((1,0),(t^\prime,v^\prime))
		\ =\  \widehat\omega_{(t,v)}(1,0)\circ\widehat\omega_{(t,v)}(t^\prime,v^\prime)
			- \widehat\omega_{(t,v)}(t^\prime,v^\prime)\circ\widehat\omega_{(t,v)}(1,0)
		\ =\ 0,$$
	since $\widehat\omega_{(t,v)}(1,0)=0$ by Remark~\ref{rem:widehatThetaOmega}.
	Therefore, we obtain
	(\ref{eqn:diffEqnWidehatOmega}).
\end{proof}
\begin{corollary}
\label{cor:diffEqnWidehatOmega}
	If the torsion tensor $\Tor$  and the curvature tensor $R$ both are parallel on $W$, i.e., $\nabla\!_v \Tor = 0$ and $\nabla\!_v R = 0$ for all vectors $v\in TW$, then we have
	\begin{eqnarray*}
		\partial_1\widehat\theta(t,v)(t^\prime,v^\prime)
		&=&
		v^\prime + \widehat\omega_{(t,v)}(t^\prime,v^\prime)(v) + \Tor_b(v,\widehat\theta_{(t,v)}(t^\prime,v^\prime)) \quad\mbox{and}\\
		\partial_1\widehat\omega(t,v)(t^\prime,v^\prime)
		&=& R_b(v,\widehat\theta_{(t,v)}(t^\prime,v^\prime),\cdot)
	\end{eqnarray*}
	for all $(t,v)\in\widehat V$ and $(t^\prime,v^\prime)\in\RR\times T_bM$.
\end{corollary}
\begin{proof}
	The proof works as in the finite-dimensional case (cf.\ \cite[p.~109]{Loo69}).
	We shall prove these assertions in two steps.
	
	{\bf Step 1:} \emph{Observe that $\Theta(v^\ast,w^\ast)(x) \equiv \Tor_b(v,w)$ and $\Omega(v^\ast,w^\ast)(x) \equiv R_b(v,w,\cdot)$ for all $x\in W$ and $v,w\in T_bM$.}
	Due to $\Theta(v^\ast,w^\ast)=\theta(\Tor(v^\ast,w^\ast))$ and $\Omega(v^\ast,w^\ast)(x)(z)=\theta(R(v^\ast,w^\ast,z^\ast))(x)$ (given $z\in T_bM$), it suffices to show that $\Tor(v^\ast,w^\ast)$ and $R(v^\ast,w^\ast,z^\ast)$ both are parallel along each geodesic $\alpha_x$ (see above) emanating from $b$, but this is true, as the curves $v^\ast\circ\alpha_x$, $w^\ast\circ\alpha_x$ and $z^\ast\circ\alpha_x$ are $\alpha_x$-parallel and as parallel transport commutes with parallel tensors (cf.\ Section~\ref{sec:affineConnections}).
	
	{\bf Step 2:} \emph{Observe that $(\Phi^\ast\Theta)_{(t,v)}((1,0),(t^\prime,v^\prime)) = \Tor_b(v,\widehat\theta_{(t,v)}(t^\prime,v^\prime))$ 	
	for all $(t,v)\in\widehat V$ and $(t^\prime,v^\prime)\in\RR\times T_bM$,
	as well as $(\Phi^\ast\Omega)_{(t,v)}((1,0),(t^\prime,v^\prime)) = R_b(v,\widehat\theta_{(t,v)}(t^\prime,v^\prime),\cdot)$.}
	By Step~1, we have
	\begin{eqnarray*}
		(\Phi^\ast\Omega)_{(t,v)}((1,0),(t^\prime,v^\prime))
		&=& \Omega_{\Phi(t,v)}(T_{(t,v)}\Phi(1,0),
					T_{(t,v)}\Phi(t^\prime,v^\prime))\\
		&=& R_b\big(\theta_{\Phi(t,v)}(T_{(t,v)}\Phi(1,0)),
		\theta_{\Phi(t,v)}(T_{(t,v)}\Phi(t^\prime,v^\prime)), \cdot\big) \\
		&\stackrel{\mbox{\scriptsize (\ref{eqn:T_tvExp(v)2})}}{=}& R_b\big(v,(\Phi^\ast\theta)_{(t,v)}(t^\prime,v^\prime), \cdot\big),
	\end{eqnarray*}
	This equals
	$R_b(v,t^\prime v + \widehat\theta_{(t,v)}(t^\prime,v^\prime),\cdot)$ by Lemma~\ref{lem:phi*thetaOmega},
	which entails the assertion, since\linebreak $t^\prime R_b(v,v,\cdot)$ vanishes, the curvature tensor being skew-symmetric with respect to the first two arguments. Similarly, we have
	\begin{eqnarray*}
		(\Phi^\ast\Theta)_{(t,v)}((1,0),(t^\prime,v^\prime))
		&=& \Theta_{\Phi(t,v)}(T_{(t,v)}\Phi(1,0),
					T_{(t,v)}\Phi(t^\prime,v^\prime))\\
		&=& \Tor_b\big(\theta_{\Phi(t,v)}(T_{(t,v)}\Phi(1,0)),
		\theta_{\Phi(t,v)}(T_{(t,v)}\Phi(t^\prime,v^\prime))\big) \\
		&\stackrel{\mbox{\scriptsize (\ref{eqn:T_tvExp(v)2})}}{=}& \Tor_b\big(v,(\Phi^\ast\theta)_{(t,v)}(t^\prime,v^\prime)\big),
	\end{eqnarray*}
	which equals $\Tor_b(v,t^\prime v + \widehat\theta_{(t,v)}(t^\prime,v^\prime))$ by Lemma~\ref{lem:phi*thetaOmega} and entails the assertion, since\linebreak $t^\prime \Tor_b(v,v)$ vanishes.
\end{proof}
%
%
%
%
\subsection{Affine Maps between Normal Neighborhoods}
\label{sec:affMaps}
Throughout this section, let $(M_1,\nabla\!_1,b_1)$ and $(M_2,\nabla\!_2,b_2)$ be affine Banach manifolds with base points. For the sake of readability, we shall usually suppress the indices for $\nabla$. Let $W_1$ and $W_2$ be normal neighborhoods of $b_1$ and $b_2$, respectively, and let $\varphi_1\colon W_1 \rightarrow V_1\subseteq T_{b_1}M =: E_1$ and $\varphi_2\colon W_2 \rightarrow V_2\subseteq T_{b_2}M =: E_2$ be the associated normal charts. We denote the maps defined in Section~\ref{sec:structureEquations} with supplementary indices. A map $f$ from a base-point containing subset of $M_1$ into $M_2$ is said to be \emph{base-point preserving} if $f(b_1)=b_2$.

The following proposition characterizes affine maps between normal neighborhoods by means of the forms $\widehat\theta$ and $\widehat\omega$.
\begin{proposition}
\label{prop:affMapsNormNeighb}
	Let $f\colon W_1 \rightarrow W_2$ be a base-point preserving smooth map.
	Then the following conditions are equivalent:
	\begin{enumerate}
		\item[\rm (a)] The map $f$ is affine.
		\item[\rm (b)] The map $f$ is compatible with adapted vector fields and covariant derivatives of adapted vector fields, i.e., we have
			\begin{eqnarray*}
				Tf \circ w^\ast &=& (T_{b_1}f(w))^\ast\circ f
					 \quad\mbox{and} \\
				T_xf(\nabla\!_v w^\ast)
					&=& \nabla\!_{T_xf(v)}(T_{b_1}f(w))^\ast
			\end{eqnarray*}
			for all $x\in W_1$, $v\in T_xM_1$ and $w \in T_{b_1}M_1$.
		\item[\rm (c)] We have
			\begin{eqnarray*}
				T_{b_1}f \circ (\theta_1)_x
					&=& (f^\ast\theta_2)_x \quad \mbox{and} \\
				T_{b_1}f \circ (\omega_1)_x(v)
					&=& (f^\ast\omega_2)_x(v)\circ T_{b_1}f
			\end{eqnarray*}
			for all $x\in W_1$ and vectors $v\in T_xM_1$.
		\item[\rm (d)] We have
			$$(T_{b_1}f)(V_1) \subseteq V_2 \quad \mbox{and}\quad
				f \circ \exp_{b_1}|{V_1} = \exp_{b_2}\circ T_{b_1}f|_{V_1},$$
			so that
			$$B\colon \widehat V_1 \rightarrow \widehat V_2,\ (t,v) \mapsto (t,T_{b_1}f(v))$$
			defines a map.
			The equations
			\begin{eqnarray*}
				 T_{b_1}f \circ (\widehat\theta_1)_{(t,v)}
				 	&=& (B^\ast\widehat\theta_2)_{(t,v)}
				 	\quad \mbox{and} \\
				 T_{b_1}f \circ (\widehat\omega_1)_{(t,v)}(t^\prime,v^\prime)
				 	&=& (B^\ast\widehat\omega_2)_{(t,v)}(t^\prime,v^\prime)\circ T_{b_1}f
			\end{eqnarray*}
			hold for all $(t,v)\in \widehat V_1$
			and $(t^\prime,v^\prime)\in \RR\times T_{b_1}M_1$.
	\end{enumerate}
\end{proposition}
\begin{proof}
	Cf.\ \cite[pp.~109-110]{Loo69} for the finite-dimensional case.
	 
	(a)$\Rightarrow$(d): In a first step, we shall deduce (b). Given any $x\in W_1$, the geodesic $f\circ\alpha_x$ joins $(f\circ\alpha_x)(0)=b_2$ with $(f\circ\alpha_x)(1)=f(x)$ and is therefore
	equal to $\alpha_{f(x)}$.
	We have
	\begin{eqnarray*}
		(Tf\circ w^\ast)(x)
		& =&  Tf(P_0^1(\alpha_{x})(w))
		\ =\ P_0^1(f\circ \alpha_x)(Tf(w)) \\
		& =& P_0^1(\alpha_{f(x)})(Tf(w))
		\ =\ (Tf(w))^\ast(f(x)),
	\end{eqnarray*}
	i.e., $f$ is compatible with adapted vector fields. Then $w^\ast$ and $(T_{b_1}f(w))^\ast$ being $f$-related, also the second equation in (b) holds. We shall now deduce (c). Similarly to the preceding considerations, we obtain
	\begin{eqnarray}
		T_{b_1}f \circ (\theta_1)_x \nonumber
			&=& T_{b_1}f \circ P_1^0(\alpha_x)
		\ =\ P_1^0(\alpha_{f(x)}) \circ T_xf \nonumber \\
		&=& (\theta_2)_{f(x)}\circ T_xf 
		\ =\ (f^\ast\theta_2)_x,	\label{eqn:affMapsNormNeighb(a)=>(d)}
	\end{eqnarray}
	which shows the first equation in (c).
	Given any $w\in T_{b_1}M$, we have
	\begin{eqnarray*}
		(T_{b_1}f \circ (\omega_1)_x(v))(w)
			&=& T_{b_1}f((\theta_1)_x(\nabla\!_v w^\ast))
			\ \stackrel{\mbox{\scriptsize (\ref{eqn:affMapsNormNeighb(a)=>(d)})}}{=} \
				(\theta_2)_{f(x)}(T_xf(\nabla\!_v w^\ast)) \\
			&\stackrel{\mbox{\scriptsize (b)}}{=}&
				(\theta_2)_{f(x)}(\nabla\!_{T_xf(v)}(T_{b_1}f(w))^\ast)
			\ = \ (\omega_2)_{f(x)}(T_xf(v))(T_{b_1}f(w)) \\
			&=& (f^\ast\omega_2)_x(v)(T_{b_1}f(w)),
	\end{eqnarray*}
	which shows the second equation in (c).
	Finally, we shall deduce (d). As, for $i=1,2$, the exponential $\exp_{b_i}^{W_i}$ of $W_i$ is given by $\exp_{b_i}|_{V_i}$, we have
	$$(T_{b_1}f)(V_1) \subseteq V_2 \quad \mbox{and}\quad
				f \circ \exp_{b_1}|{V_1} = \exp_{b_2}\circ T_{b_1}f|_{V_1}.$$
	A simple computation shows $\Phi_2\circ B = f\circ\Phi_1$.
	Hence, it follows that
	\begin{eqnarray*}
		T_{b_1f} \circ (\Phi_1^\ast\theta_1)_{(t,v)}
			&=& T_{b_1}f \circ (\theta_1)_{\Phi_1(t,v)}\circ T_{(t,v)}\Phi_1
		\ \stackrel{\mbox{\scriptsize (c)}}{=} \
			(f^\ast\theta_2)_{\Phi_1(t,v)}\circ T_{(t,v)}\Phi_1 \\
		&=& (\Phi_1^\ast(f^\ast\theta_2))_{(t,v)}
		\ =\ (B^\ast(\Phi_2^\ast\theta_2))_{(t,v)},
	\end{eqnarray*}
	i.e.,
	$$T_{b_1f} \circ (\lambda_{T_{b_1}M_1}\,d\lambda_{\RR} + \widehat\theta_1)_{(t,v)} \ =\ (B^\ast(\lambda_{T_{b_2}M_2}\,d\lambda_{\RR} + \widehat\theta_2))_{(t,v)},$$
	by Lemma~\ref{lem:phi*thetaOmega}.
	As it is easy to check that
	\begin{equation} \label{eqn:affMapsNormNeighb(a)=>(d)2}
		T_{b_1}f \circ (\lambda_{T_{b_1}M_1}\,d\lambda_{\RR})_{(t,v)}
	\ =\
	(B^\ast(\lambda_{T_{b_2}M_2}\,d\lambda_{\RR}))_{(t,v)},
	\end{equation}
	we obtain the first equation in (d).
	The second one follows by
	\begin{eqnarray*}
		T_{b_1}f\circ(\Phi_1^\ast\omega_1)_{(t,v)}(t^\prime,v^\prime)
			&=& T_{b_1}f\circ
				(\omega_1)_{\Phi_1(t,v)}(T_{(t,v)}\Phi_1(t^\prime,v^\prime)) \\
		& \stackrel{\mbox{\scriptsize (c)}}{=} &
				(f^\ast\omega_2)_{\Phi_1(t,v)}(T_{(t,v)}\Phi_1(t^\prime,v^\prime))
				\circ T_{b_1}f \\
		&=& (\Phi_1^\ast(f^\ast\omega_2))_{(t,v)}(t^\prime,v^\prime) \circ T_{b_1}f
		\ = \ 
		(B^\ast(\Phi_2^\ast\omega_2))_{(t,v)}(t^\prime,v^\prime) \circ T_{b_1}f,
	\end{eqnarray*}
	since we have $\Phi_1^\ast\omega_1 = \widehat\omega_1$ and $\Phi_2^\ast\omega_2 = \widehat\omega_2$ by Lemma~\ref{lem:phi*thetaOmega}.

	(d)$\Rightarrow$(c): By adding  (\ref{eqn:affMapsNormNeighb(a)=>(d)2})
	and $T_{b_1}f \circ (\widehat\theta_1)_{(t,v)} = (B^\ast\widehat\theta_2)_{(t,v)}$, we get
	$$T_{b_1}f \circ(\Phi_1^\ast\theta_1)_{(t,v)}
		\ =\ (B^\ast(\Phi_2^\ast\theta_2))_{(t,v)},$$
	for all $(t,v)\in \widehat V_1$, due to $\lambda_{T_{b_1}M_1}\,d\lambda_{\RR} + \widehat\theta_1 = \Phi_1^\ast\theta_1$ and $\lambda_{T_{b_2}M_2}\,d\lambda_{\RR} + \widehat\theta_2  = \Phi_2^\ast\theta_2$ (cf.\ Lemma~\ref{lem:phi*thetaOmega}).
	As above, we have $\Phi_2\circ B = f\circ\Phi_1$, so that $B^\ast(\Phi_2^\ast\theta_2)=\Phi_1^\ast(f^\ast\theta_2)$ and hence
	$$T_{b_1}f \circ (\theta_1)_{\Phi_1(t,v)}\circ T_{(t,v)}\Phi_1
		\ =\ (f^\ast\theta_2)_{\Phi_1(t,v)} \circ T_{(t,v)}\Phi_1$$
	for all $(t,v) \in \widehat V_1$.
	In order to show
	$T_{b_1}f \circ (\theta_1)_x = (f^\ast\theta_2)_x,$
	we put $t:=1$ and $v:=\varphi_1(x)$, which entails $\Phi_1(t,v)= x$.
	Then we have
	$$T_{b_1}f \circ (\theta_1)_x\circ T_{(1,v)}\Phi_1
		\ =\ (f^\ast\theta_2)_x \circ T_{(1,v)}\Phi_1.$$
	Therefore, it suffices to show the surjectivity of $T_{(1,v)}\Phi_1$. In fact, even the partial map $T_{(1,v)}\Phi_1(0,\cdot)\stackrel{\mbox{\scriptsize (\ref{eqn:TPhi})}}{=}T_v\exp_{b_1}$ is surjective, the exponential $\exp_{b_1}|_{V}^{W}$ being a diffeomorphism.
	Let us now deduce the second equation in (c).
	From the second equation
	in (d) and from $\Phi_1^\ast\omega_1 = \widehat\omega_1$ and $\Phi_2^\ast\omega_2 = \widehat\omega_2$ (cf.\ Lemma~\ref{lem:phi*thetaOmega}), we immediately obtain
	$$T_{b_1}f \circ (\Phi_1^\ast\omega_1)_{(t,v)}(t^\prime,v^\prime)
		\ =\ (B^\ast(\Phi_2^\ast\omega_2))_{(t,v)}(t^\prime,v^\prime)\circ T_{b_1}f$$
	for all $(t,v)\in\widehat V$ and $(t^\prime,v^\prime)\in \RR\times T_{b_1}M_1$.
	Due to $B^\ast(\Phi_2^\ast\omega_2)=\Phi_1^\ast(f^\ast\omega_2)$, we then have
	$$T_{b_1}f \circ (\omega_1)_{\Phi_1(t,v)}
			(T_{(t,v)}\Phi_1(t^\prime,v^\prime))
		\ =\ (f^\ast\omega_2)_{\Phi_1(t,v)}(T_{(t,v)}\Phi_1(t^\prime,v^\prime))\circ T_{b_1}f.$$
	In order to show
	\begin{equation*} 
		T_{b_1}f \circ (\omega_1)_x(w) \ =\ (f^\ast\omega_2)_x(w)\circ T_{b_1}f
	\end{equation*}
	for arbitrary $x\in W_1$ and $w\in T_xM_1$,
	we put $t:=1$ and $v:=\varphi_1(x)$, which entails $\Phi_1(t,v)=x$. Further, we put $t^\prime:=0$ and 
	choose $v^\prime$ in such a way that $T_{(1,v)}\Phi_1(0,v^\prime)=w$, which is possible, the map $T_{(1,v)}\Phi_1(0,\cdot)$ being an isomorphism (see above).

	(c)$\Rightarrow$(b): Given any $x\in W_1$, we deduce from the first equation in (c) that
	$T_{b_1}f \circ P_1^0(\alpha_x) = P_1^0(\alpha_{f(x)}) \circ T_xf,$
	which is equivalent to
	$P_0^1(\alpha_{f(x)}) \circ T_{b_1}f = T_xf \circ P_0^1(\alpha_x)$.
	Evaluating the two sides at any $w\in T_{b_1}M_1$ and respecting that this holds for all $x\in W_1$, we obtain
	$(T_{b_1}f(w))^\ast\circ f = Tf\circ w^\ast$,
	i.e., the first equation in (b). We shall deduce the second one.
	Given any $x\in W_1$, $v\in T_xM_1$ and $w\in T_{b_1}M_1$, we evaluate both sides of
	$T_{b_1}f \circ (\omega_1)_x(v)
		= (f^\ast\omega_2)_x(v)\circ T_{b_1}f$
	at $w$ and obtain
	$$T_{b_1}f((\theta_1)_x(\nabla\!_v w^\ast))
		\ =\ (\theta_2)_{f(x)}(\nabla\!_{T_xf(v)}(T_{b_1}f(w))^\ast).$$
	Since
	$T_{b_1}f \circ (\theta_1)_x = (f^\ast\theta_2)_x$ by (c), the left hand side equals
	$(\theta_2)_{f(x)}(T_xf(\nabla\!_v w^\ast))$, so that we have
	$$T_xf(\nabla\!_v w^\ast)
		\ =\ \nabla\!_{T_xf(v)}(T_{b_1}f(w))^\ast,$$
	the map $(\theta_2)_{f(x)}$ being an isomorphism.

	(b)$\Rightarrow$(a): Given any pair $(\varphi_1,\varphi_2)$ of charts $\varphi_1\colon U_{\varphi_1} \rightarrow V_{\varphi_1}\subseteq E_1$ of $W_1$ and\linebreak $\varphi_2\colon U_{\varphi_2} \rightarrow V_{\varphi_2}\subseteq E_2$ of $W_2$  such that $f(U_{\varphi_1})\subseteq U_{\varphi_2}$, we shall verify that
	\begin{equation} \label{eqn:affMapsNormNeighb(b)=>(a)}
	d^2f^\varphi(\bar x)(\bar v,\bar w)
		+ df^\varphi(\bar x)((B_1^{\varphi_1})_{\bar x}
			(\bar v,\bar w))
		\ =\ (B_2^{\varphi_2})_{f^\varphi(\bar x)}
				(df^\varphi(\bar x)(\bar v),
				df^\varphi(\bar x)(\bar w))
	\end{equation}
	for all $\bar x\in V_{\varphi_1}$ and $\bar v,\bar w\in E_1$.
	We put $x:=\varphi_1^{-1}(\bar x)$, $v:=T\varphi_1^{-1}(\bar x,\bar v)$ and $w:=\linebreak T\varphi_1^{-1}(\bar x,\bar w)$.
	From (b), we know that the smooth vector fields
	$\eta_1:= ((\theta_1)_x(v))^\ast$ and $\eta_2:= \big(T_{b_1}f((\theta_1)_x(v))\big)^\ast$
	satisfy
	$$Tf \circ \eta_1 =  \eta_2 \circ f \quad\mbox{and}\quad T_xf(\nabla\!_w\eta_1) = \nabla\!_{T_xf(w)}\eta_2,$$
	so that working in the charts lets us know
	$$d^2f^\varphi(\bar x)(\eta_1^{\varphi_1}(\bar x),\bar w)
			+ df^\varphi(\bar x)(d\eta_1^{\varphi_1}(\bar x)(\bar w))
		\ =\ d\eta_2^{\varphi_2}(f^\varphi(\bar x))(df^\varphi(\bar x)(\bar w))$$
	(where we have taken the derivative at $\bar x$ in direction $\bar w$) and
	\begin{eqnarray*}
		\lefteqn{df^\varphi(\bar x)\big(d\eta_1^{\varphi_1}(\bar x)(\bar w)
		- (B_1^{\varphi_1})_{\bar x}(\eta_1^{\varphi_1}(\bar x),\bar w)\big)}  \\
		&\quad=& d\eta_2^{\varphi_2}(f^\varphi(\bar x))
				(df^\varphi(\bar x)(\bar w))
			- (B_2^{\varphi_2})_{f^\varphi(\bar x)}
				\big(\eta_2^{\varphi_2}(f^\varphi(\bar x)),
				df^\varphi(\bar x)(\bar w)\big).
	\end{eqnarray*}
	By subtracting the latter equation from the former one, we obtain
	\begin{eqnarray*}
		\lefteqn{d^2f^\varphi(\bar x)(\eta_1^{\varphi_1}(\bar x),\bar w)
			+ df^\varphi(\bar x)\big((B_1^{\varphi_1})_{\bar x}
				(\eta_1^{\varphi_1}(\bar x),\bar w)\big)} \\
			&\quad =& (B_2^{\varphi_2})_{f^\varphi(\bar x)}
				\big(\eta_2^{\varphi_2}(f^\varphi(\bar x)),
				df^\varphi(\bar x)(\bar w)\big),
	\end{eqnarray*}
	which leads to (\ref{eqn:affMapsNormNeighb(b)=>(a)}), since
	$\eta_1^{\varphi_1}(\bar x) = \bar v$ and $\eta_2^{\varphi_2}(f^\varphi(\bar x)) = df^\varphi(\bar x)(\bar v)$.
\end{proof}
\begin{theorem}
\label{th:affMapsInfCharac}
	Assume that the torsion and curvature tensors $\Tor_1$, $\Tor_2$, $R_1$ and $R_2$ of $(M_1,\nabla\!_1,b_1)$ and $(M_2,\nabla\!_2,b_2)$ are parallel on $W_1$ and $W_2$, respectively, i.e., we have $\nabla\!_v\Tor_1=0$ and $\nabla\!_vR_1=0$ for all $v\in TW_1$ and $\nabla\!_v\Tor_2=0$ and $\nabla\!_vR_2=0$ for all $v\in TW_2$.
	Then, for each base-point preserving smooth map $f\colon W_1 \rightarrow W_2$, the following conditions are equivalent:
	\begin{enumerate}
		\item[\rm (a)] The map $f$ is affine.
		\item[\rm (b)] The map $f$ intertwines the torsion and curvature tensors in the base points and the exponential maps $\exp_{b_1}|_{V_1}$ and $\exp_{b_2}|_{V_2}$, i.e., we have
			$$T_{b_1}f \circ (\Tor_1)_{b_1} = (\Tor_2)_{b_2}\circ (T_{b_1}f)^2 \quad\mbox{and}\quad T_{b_1}f \circ (R_1)_{b_1} = (R_2)_{b_2}\circ (T_{b_1}f)^3$$
			as well as
			$$(T_{b_1}f)(V_1) \subseteq V_2 \quad \mbox{and}\quad
				f \circ \exp_{b_1}|{V_1} = \exp_{b_2}\circ T_{b_1}f|_{V_1}.$$
	\end{enumerate}
\end{theorem}
\begin{proof}
	The proof works as in the finite-dimensional case (cf.\ \cite[p.~111]{Loo69}). Note that O.~Loos there claims a stronger result (cf.\ Rem.~\ref{rem:LoosGap}).
		
	(a)$\Rightarrow$(b) is clear.
	
	(b)$\Rightarrow$(a): It suffices to show Condition (d) of Proposition~\ref{prop:affMapsNormNeighb}, i.e.,
	\begin{eqnarray*}
		T_{b_1}f \cdot \widehat\theta_1
			&=& B^\ast\widehat\theta_2
				\quad \mbox{and} \\
		T_{b_1}f \cdot \widehat\omega_1
			&=& B^\ast\widehat\omega_2 \cdot T_{b_1}f,
	\end{eqnarray*}
	where $\cdot$ denotes the appropriate maps.
	To verify these equations, we pursue the idea to show that the pairs
	$(T_{b_1}f \cdot \widehat\theta_1, T_{b_1}f \cdot \widehat\omega_1)$
	and $(B^\ast\widehat\theta_2,B^\ast\widehat\omega_2 \cdot T_{b_1}f)$
	satisfy the same system of ordinary differential equations with initial conditions. Then they are equal by the Uniqueness Theorem (cf.\ \cite[pp.~70, 72]{Lan01}).
	We claim that $(T_{b_1}f \cdot \widehat\theta_1, T_{b_1}f \cdot \widehat\omega_1)$ satisfies the two equations
	\begin{eqnarray}
		\lefteqn{\partial_1 (T_{b_1}f \cdot \widehat\theta_1)(t,v)(t^\prime,v^\prime)} \label{eqn:affMapsInfCharac1}\\
			&\quad = &  T_{b_1}f(v^\prime)
				+ (T_{b_1}f \cdot \widehat\omega_1)(t,v)(t^\prime,v^\prime)(v)
				+ (\Tor_2)_{b_2}(T_{b_1}f(v),(T_{b_1}f\cdot\widehat\theta_1)(t,v)(t^\prime,v^\prime))	 \nonumber\\
		\lefteqn{\partial_1(T_{b_1}f \cdot \widehat\omega_1)(t,v)(t^\prime,v^\prime)
			\ =\ (R_2)_{b_2}(T_{b_1}f(v),
					(T_{b_1}f \cdot \widehat\theta_1)(t,v)(t^\prime,v^\prime),
					\cdot) \circ T_{b_1}f}
				\label{eqn:affMapsInfCharac2}
	\end{eqnarray}
	with initial conditions
	$(T_{b_1}f \cdot \widehat\theta_1)(0,v)=0$ and
	$(T_{b_1}f \cdot \widehat\omega_1)(0,v) = 0$
	for all $(t,v)\in\widehat V_1$ and $(t^\prime,v^\prime)\in \RR\times T_{b_1}M_1$.
	We further claim that $(B^\ast\widehat\theta_2,B^\ast\widehat\omega_2 \cdot T_{b_1}f)$ satisfies the two equations
	\begin{eqnarray}
		\lefteqn{\partial_1 (B^\ast\widehat\theta_2)(t,v)(t^\prime,v^\prime)} \label{eqn:affMapsInfCharac4} \\
			&\quad =&  T_{b_1}f(v^\prime)
				+ (B^\ast\widehat\omega_2\cdot T_{b_1}f)(t,v)(t^\prime,v^\prime)(v)
				+ (\Tor_2)_{b_2}(T_{b_1}f(v),(B^\ast\widehat\theta_2)(t,v)(t^\prime,v^\prime)) \nonumber \\
		\lefteqn{\partial_1(B^\ast\widehat\omega_2 \cdot
			T_{b_1}f)(t,v)(t^\prime,v^\prime)
			\ =\ (R_2)_{b_2}(T_{b_1}f(v),
				(B^\ast\widehat\theta_2)(t,v)(t^\prime,v^\prime),\cdot) \circ T_{b_1}f}
					\label{eqn:affMapsInfCharac5}
	\end{eqnarray}
	with initial conditions
	$B^\ast\widehat\theta_2(0,v)=0$ and
	$(B^\ast\widehat\omega_2\cdot T_{b_1}f)(0,v) = 0$.
	It suffices to check the equations $(\ref{eqn:affMapsInfCharac1})$ - $(\ref{eqn:affMapsInfCharac5})$. Note that the initial conditions follow by the mere definitions of the respective maps.
	
	By Corollary~\ref{cor:diffEqnWidehatOmega}, we have
	\begin{eqnarray*}
		\partial_1 (T_{b_1}f \cdot \widehat\theta_1)(t,v)(t^\prime,v^\prime)
		&=& (T_{b_1}f \circ \partial\widehat\theta_1(t,v))(t^\prime,v^\prime) \\
		&=& T_{b_1}f\big(v^\prime + (\widehat\omega_1)_{(t,v)}(t^\prime,v^\prime)(v) + (\Tor_1)_{b_1}(v,(\widehat\theta_1)_{(t,v)}(t^\prime,v^\prime))\big),
	\end{eqnarray*}
	which equals
	$$T_{b_1}f\big(v^\prime) + (T_{b_1}f\cdot\widehat\omega_1)(t,v)(t^\prime,v^\prime)(v) + (\Tor_2)_{b_2}(T_{b_1}f(v),(T_{b_1}f\cdot\widehat\theta_1)(t,v)(t^\prime,v^\prime))\big),$$
	as the map $f$ intertwines the torsion tensors in the base points. Thus (\ref{eqn:affMapsInfCharac1}) is shown.
	Similarly, (\ref{eqn:affMapsInfCharac2}) holds, since we have
	\begin{eqnarray*}
		\partial_1 (T_{b_1}f \cdot \widehat\omega_1)(t,v)(t^\prime,v^\prime)
		&=& T_{b_1}f \circ \partial\widehat\omega_1(t,v)(t^\prime,v^\prime) \\
		&=& T_{b_1}f\circ
			(R_1)_{b_1}(v,(\widehat\theta_1)_{(t,v)}(t^\prime,v^\prime),\cdot) \\
		&=& (R_2)_{b_2}(T_{b_1}f(v),
					 (T_{b_1}f\cdot\widehat\theta_1)(t,v)(t^\prime,v^\prime),
					\cdot) \circ T_{b_1}f,
	\end{eqnarray*}
	the map $f$ intertwining the curvature tensors in the base points.
		
	The map $B$ being a restriction of the continuous linear map $\id_{\RR}\times T_{b_1}f$, we have
	$$\partial_1(B^\ast\widehat\theta_2)(t,v)(t^\prime,v^\prime)
		\ =\ (B^\ast(\partial_1\widehat\theta_2))(t,v)(t^\prime,v^\prime)
		\ =\ \partial_1\widehat\theta_2(B(t,v))(t^\prime, T_{b_1}f(v^\prime)),$$
	which equals
	$$T_{b_1}f(v^\prime) + (\widehat\omega_2)_{B(t,v)}(t^\prime,T_{b_1}f(v^\prime))(T_{b_1}f(v)) + (\Tor_2)_{b_2}\big(T_{b_1}f(v),(\widehat\theta_2)_{B(t,v)}(t^\prime,T_{b_1}f(v^\prime))\big),$$
	by Corollary~\ref{cor:diffEqnWidehatOmega}, i.e.,
	$$T_{b_1}f(v^\prime)
		+ (B^\ast\widehat\omega_2\cdot T_{b_1}f)(t,v)(t^\prime,v^\prime)(v)
		+ (\Tor_2)_{b_2}(T_{b_1}f(v),(B^\ast\widehat\theta_2)(t,v)(t^\prime,v^\prime)).$$
	Thus (\ref{eqn:affMapsInfCharac4}) is shown.
	Similarly, (\ref{eqn:affMapsInfCharac5}) holds, since we have
	\begin{eqnarray*}
		\partial_1(B^\ast\widehat\omega_2 \cdot T_{b_1}f)(t,v)(t^\prime,v^\prime)
		& =& (B^\ast(\partial_1\widehat\omega_2))(t,v)(t^\prime,v^\prime)\circ T_{b_1}f \\
		&=& \partial_1\widehat\omega_2(B(t,v))(t^\prime, T_{b_1}f(v^\prime))\circ T_{b_1}f \\
		&=& (R_2)_{b_2}\big(T_{b_1}f(v),(\widehat\theta_2)_{B(t,v)}(t^\prime,T_{b_1}f(v^\prime)),\cdot\big)\circ T_{b_1}f \\
		&=& (R_2)_{b_2}(T_{b_1}f(v),
				(B^\ast\widehat\theta_2)(t,v)(t^\prime,v^\prime),\cdot) \circ T_{b_1}f.
	\end{eqnarray*}
\end{proof}
\begin{corollary}[Local Integrability] \label{cor:extensionToAffMaps}
	We assume that the torsion and curvature tensors $\Tor_1$, $\Tor_2$, $R_1$ and $R_2$ of $(M_1,\nabla\!_1,b_1)$ and $(M_2,\nabla\!_2,b_2)$ are parallel on certain neighborhoods of $b_1$ and $b_2$, respectively.
	For every continuous linear map $A\colon T_{b_1}M_1 \rightarrow T_{b_2}M_2$ that intertwines the torsion and curvature tensors in the base points in the sense that
	$$A\circ (\Tor_1)_{b_1} = (\Tor_2)_{b_2}\circ A^2 \quad\mbox{and}\quad A \circ (R_1)_{b_1} = (R_2)_{b_2}\circ A^3,$$
	there exists an affine map $f$ from an open neighborhood of $b_1$ into $M_2$ that satisfies $T_{b_1}f = A$.
\end{corollary}
\begin{proof}
	We make the normal charts $\varphi_1\colon W_1\rightarrow V_1\subseteq E_1$ and $\varphi_2\colon W_2\rightarrow V_2\subseteq E_2$ sufficiently small so that $R_1$ and $R_2$ are parallel on $W_1$ and $W_2$, respectively. If necessary, we again shrink $W_1$ such that
	$A(V_1) \subseteq V_2$.
	We define 
	$f\colon W_1 \rightarrow W_2$
	by $f:=\exp_{b_2}\circ A \circ \varphi_1$
	and have
	\begin{equation} \label{eqn:extensionToAffMaps}
	f\circ\exp_{b_1}|_{V_1}
		\ =\ \exp_{b_2}\circ A|_{V_1},
	\end{equation}
	due to $\varphi_1=(\exp_{b_1}|_{V_1}^{W_1})^{-1}$.
	It suffices to check $T_{b_1}f=A$, as then the map $f$ is affine by Theorem~\ref{th:affMapsInfCharac}.
	By taking the derivative of both sides of (\ref{eqn:extensionToAffMaps}) at $0_{b_1}$ in any direction $v\in T_{b_1}M_1$, we obtain
	$$T_{b_1}f(T_{0_{b_1}}\!\exp_{b_1}(v))
		\ =\ T_{A(0_{b_1})}\exp_{b_2}(dA(0_{b_1})(v)).$$
	Due to $T_{0_{b_i}}\!\exp_{b_i}=\id_{T_{b_i}M_i}$ for $i=1,2$, we then have
	$T_{b_1}f(v)= A(v)$,
	the map $A$ being continuous linear. Thus, $T_{b_1}f=A$ indeed holds.
\end{proof}
\begin{remark}\label{rem:LoosGap}
	In \cite{Loo69}, O.~Loos claims that in the finite-dimensional case, every linear map $A$ that intertwines the torsion and curvature tensors in the base points can be extended to an affine map on any given normal neighborhood provided $(M_2,\nabla_2,b_2)$ is complete. Regrettably, his argument seems to be incomplete.
\end{remark}
%
%
%
%
%
%
%
\section{Integration of Maps between Tangent Spaces to Affine Maps}
\label{sec:integrationToAffMaps}
Let $(M_1,\nabla\!_1,b_1)$ and $(M_2,\nabla\!_2,b_2)$ be affine Banach manifolds with base points. A continuous linear map $A\colon T_{b_1}M_1\rightarrow T_{b_2}M_2$ between the tangent spaces at the base points is called \emph{integrable} if there exists an affine map $f\colon M_1 \rightarrow M_2$ that satisfies $T_{b_1}f=A$. If $M_1$ is connected, then the map $f$ is unique if it exists.

When knowing that $A$ is integrable on a neighborhood of $b_1$, we aim at global integrability by the extension of local integrals along piecewise geodesics. The main result of this section is that for $1$-connected $M_1$ and geodesically complete $M_2$, a locally integrable map\linebreak $A\colon T_{b_1}M_1\rightarrow T_{b_2}M_2$ is integrable if and only if, in plain terms, extension along piecewise geodesics is possible. 

Combining this result with the theorem about local integrability, we observe that a continuous linear map $A\colon T_{b_1}M_1\rightarrow T_{b_2}M_2$ that intertwines the torsion and curvature tensors in the base points is integrable if the manifolds have parallel torsion and curvature and if $M_1$ is $1$-connected and $M_2$ is geodesically complete.

\subsection{Piecewise Geodesics}
\label{sec:piecewiseGeodesics}
Let $(M,\nabla)$ be an affine Banach manifold. We recall that a geodesic in $M$ is essentially determined by a single point and its respective velocity vector. Different velocity vectors along the geodesic are related by parallel transport.
 
A \emph{piecewise geodesic in $M$} is a continuous curve $\alpha\colon I=[a,b]\rightarrow M$ for which there is a partition of $I$ into intervals $I_1=[t_0,t_1]$, $I_2=[t_1,t_2]$, \dots, $I_n=[t_{n-1},t_n]$ such that the restrictions $\alpha_1:=\alpha|_{I_1}$, \dots, $\alpha_n:=\alpha|_{I_n}$ are geodesics in $M$. We shall often use the notation $\alpha=(\alpha_1,\dots,\alpha_n)$.

The curve can be described by the following data:
by a distinguished time $d\in I$ and its respective point $\alpha(d)$, by the time points $t_0,\dots,t_n$ and by the vectors $v_1,\dots,v_n\in T_{\alpha(d)}M$ that satisfy
$$\alpha_i^\prime(t)=P_{d}^{t}(\alpha)(v_i)\quad\mbox{for all $t \in I_i=[t_{i-1},t_i]$ with $i=1,\dots, n$},$$
i.e., by the data $(d,\alpha(d);t_0\leq t_1\leq\ldots\leq t_n;v_1,\dots,v_n)$. The curve is uniquely determined by this data.

Conversely, given data $(d,x;t_0\leq t_1\leq\ldots\leq t_n;v_1,\dots,v_n)$ with $d\in [t_0,t_n]$, $x\in M$ and $v_1,\dots, v_n\in T_xM$, there need not be a piecewise geodesic that can be described by them, but there is one if $(M,\nabla)$ is geodesically complete.

The following lemma is not surprising:
\begin{lemma}[Change of data]\label{lem:changeOfData}
	Given a piecewise geodesic $\alpha$ that can be described by the data $(d,\alpha(d);t_0\leq t_1\leq\ldots\leq t_n;v_1,\dots,v_n)$, then for each $d^\prime\in [t_0,t_n]$, it can also be described by the data $(d^\prime,\alpha(d^\prime);t_0\leq t_1\leq\ldots\leq t_n;P_{d}^{d^\prime}(\alpha)(v_1),\dots,P_{d}^{d^\prime}(\alpha)(v_n))$.
\end{lemma}
\begin{lemma}[Geodesic connection] \label{lem:geodesicConnection}
	If $M$ is connected, then, given any points $x$ and $y$ in $M$, there is a piecewise geodesic $\alpha\colon [0,1] \rightarrow M$ that joins $\alpha(0)=x$ with $\alpha(1)=y$.
\end{lemma}
\begin{proof}
	Due to the possibility of reparametrization, we need not mind the domains of the considered piecewise geodesics. Given some $x\in M$, let $A$ be the set of all $y\in M$ for which there is a piecewise geodesic that joins $x$ and $y$. To see that $A$ is all of $M$, we shall show that $A$ and its complement $A^c$ both are open. For each $y\in M$, we can consider a normal neighborhood $U_y$ and recall that $y$ can be joined with each of its points by a geodesic. Hence, if $y\in A$ then $U_y\subseteq A$, and if $y\in A^c$ then $U_y\subseteq A^c$. This shows that $A$ and $A^c$ both are open.
\end{proof}
The following lemma says that every point in $M$ has a neighborhood in which points can be joined by geodesics (that lie in $M$) such that the geodesics depend continuously on the points. Moreover, geodesics that lie in the neighborhood are uniquely determined by their endpoints.
\begin{lemma}\label{lem:UPhiW}
	Let $\Phi\colon TM\supseteq \calD_{\exp} \rightarrow M\times M$ be the map defined by
	$\Phi(v):=(\pi(v),\exp(v))$, where $\pi\colon TM\rightarrow M$ denotes the natural projection of the tangent bundle. For each $x\in M$, we have:
	\begin{enumerate}
		\item[\rm (1)]
		There exists a triple $(U_x,\Phi_x,W_x)$, where $U_x$ is an open neighborhood of $0_x$ in $\calD_{\exp}\subseteq TM$ that contains the zero section of $\pi(U_x)$,  $\Phi_x:=\Phi|_{U_x}$ is a diffeomorphism onto its open image and $W_x$ is an open neighborhood of $x$ such that $\Phi(U_x)\supseteq W_x\times W_x$.
		\item[\rm (2)]
		Given any geodesic $\alpha\colon[0,1]\rightarrow M$ that lies in $W_x$ and joins $\alpha(0)=y$ with $\alpha(1)=z$, we have $\alpha(t)=\exp_y(t\Phi_x^{-1}(y,z))$ for all $t\in [0,1]$.
	\end{enumerate}
\end{lemma}
\begin{proof}
	(1) is an easy consequence of the fact that $\Phi$ induces a local diffeomorphism at $0_x$, since $T_{0_x}\Phi$ is a topological linear isomorphism (cf.\ \cite[Prop.~5.1]{Lan01}).
	
	(2) Let $v$ be in $T_yM$ such that $\alpha(t)=\exp_y(tv)$. By definition, we have $\Phi(v)=(y,z)$. It suffices to show that $v\in U_x$, as then $v=\Phi_x^{-1}(y,z)$. For this, let $A$ be the set of all $t\in [0,1]$ for which $tv\in U_x$. We shall show that $A$ is all of $[0,1]$ and shall do this by checking that $A$ is open, closed and not empty. The latter is true, since $0v=0_y\in U_x$. The openness of $A$ is clear by the openness of $U_x$.
	
	To see the closedness of $A$, consider any sequence $(t_n)_{n\in\NN}$ in $A$ with $\lim_{n\rightarrow\infty} t_n=t\in [0,1]$ and check that $t\in A$. For this, we show that $tv\in U_x$, i.e., $\lim_{n\rightarrow\infty} t_nv \in U_x$. The sequence $(t_nv)_{n\in\NN}$ in $U_x$ converges in $U_x$ if and only if the sequence $(\Phi_x(t_nv))_{n\in\NN}=(y,\exp_y(t_nv))_{n\in\NN}$ converges in $\im(\Phi_x)$, but the latter is true, since $\lim_{n\rightarrow\infty}\exp_y(t_nv)=\exp_y(tv)=\alpha(t)\in W_x$.
\end{proof}
\begin{proposition}
	Let $\alpha_0,\alpha_1\colon [0,1]\rightarrow M$ be two piecewise geodesics with $\alpha_0(0)=\alpha_1(0)$ and $\alpha_0(1)=\alpha_1(1)$. If there is a homotopy $H\colon [0,1]\times [0,1]\rightarrow M$ between $H(\cdot,0)=\alpha_0$ and $H(\cdot,1)=\alpha_1$, then there exists a homotopy $H^\prime$ between them such that each curve $\alpha_s:=H^\prime(\cdot,s)$ is a piecewise geodesic.
\end{proposition}
\begin{proof}
	Given a collection $(U_x,\Phi_x,W_x)_{x\in M}$ of triples as in Lemma~\ref{lem:UPhiW}, we consider the covering $(H^{-1}(W_x))_{x\in M}$ of $[0,1]\times [0,1]$. Let $\lambda>0$ be a Lebesgue number of this covering, i.e., every subset of $[0,1]\times [0,1]$ whose diameter is less than $\lambda$ is contained in at least one of the open sets $H^{-1}(W_x)$ (cf.\ e.g.\ \cite[I.7.4]{Sch75}). We decompose $[0,1]$ into $0=t_0<t_1<\ldots<t_n=1$ such that, for all $i=1,\ldots,n$, the restrictions of $\alpha_0$ and $\alpha_1$ to $[t_{i-1},t_i]$ are geodesics and $|t_{i-1}-t_i|<\frac{\lambda}{2}$. 
	
	We define $H^\prime$ on $A:=\big(\{t_0,\ldots,t_n\}\times [0,1]\big) \cup \big([0,1]\times \{0,1\}\big)$ by $H^\prime|_A:=H|_A$ and target a suitable extension. On each filled rectangle $\overline R_i:=[t_{i-1},t_i]\times [0,1]$, we need a continuous extension of $H^\prime|_A$ such that $(H^\prime|_{\overline R_i})(\cdot,s)$ is a piecewise geodesic for each $s\in [0,1]$. Without loss of generality, we shall inspect $\overline R_1=[0,t_1]\times [0,1]$. For this, we consider an equidistant decomposition $0=s_0 <s_1<\ldots<s_m$ of $[0,1]$ for which each parallelogram $P_j:=P\big((0,s_j),(t_1,s_{j-1}),(t_1,s_j),(0,s_{j+1})\big)$ with $j=1,\ldots,m-1$ has diameter less than $\lambda$. The triangles $\Delta_0:=\Delta\big((0,0),(t_1,0),(0,s_1)\big)$ and $\Delta_1:=\Delta\big((0,1),(t_1,s_{m-1}),(t_1,1)\big)$ then both have diameters less than $\lambda$, too. Let $B$ be the union of all these parallelograms and triangles and note that $\overline R_1$ is their union when considering the \emph{filled} polygons. We extend $H^\prime|_A$ to $B$ by $H^\prime|_B:=H|_B$. On each filled parallelogram and triangle, we need a suitable continuous extension.
	
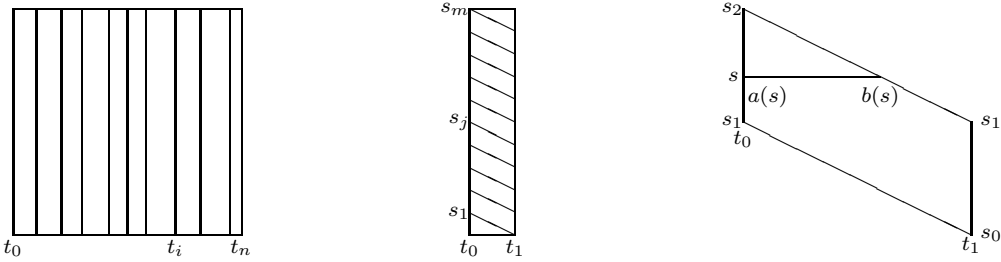
\begin{figure}[here]
\setlength{\unitlength}{3cm}
\begin{center}
\begin{picture}(4.2,1)
	\put(0,0){\line(1,0){1}}
	\put(0,1){\line(1,0){1}}
	\put(0,0){\line(0,1){1}}
	\put(0.1,0){\line(0,1){1}}
	\put(0.21,0){\line(0,1){1}}
	\put(0.3,0){\line(0,1){1}}
	\put(0.42,0){\line(0,1){1}}
	\put(0.5,0){\line(0,1){1}}
	\put(0.58,0){\line(0,1){1}}
	\put(0.71,0){\line(0,1){1}}
	\put(0.82,0){\line(0,1){1}}
	\put(0.95,0){\line(0,1){1}}
	\put(1,0){\line(0,1){1}}
	%
	\put(0,-0.02){\makebox(0,0)[t]{\scriptsize $t_0$}}
	\put(0.71,-0.02){\makebox(0,0)[t]{\scriptsize $t_i$}}
	\put(1,-0.02){\makebox(0,0)[t]{\scriptsize $t_n$}}
	%
	%
	\put(2,0){\line(1,0){0.2}}
	\put(2,1){\line(1,0){0.2}}
	\put(2,0){\line(0,1){1}}
	\put(2.2,0){\line(0,1){1}}
	\multiput(2,0.1)(0,0.1){10}{\line(2,-1){0.2}}
	\put(2,-0.02){\makebox(0,0)[t]{\scriptsize $t_0$}}
	\put(2.2,-0.02){\makebox(0,0)[t]{\scriptsize $t_1$}}
	\put(2,0.1){\makebox(0,0)[r]{\scriptsize $s_1$}}
	\put(2,0.5){\makebox(0,0)[r]{\scriptsize $s_j$}}
	\put(2,1){\makebox(0,0)[r]{\scriptsize $s_m$}}
	%
	%
	\put(3.2,0.5){\line(0,1){0.5}}
	\put(4.2,0){\line(0,1){0.5}}
	\put(3.2,0.5){\line(2,-1){1}}
	\put(3.2,1){\line(2,-1){1}}
	\put(3.2,0.7){\line(1,0){0.6}}
	\put(3.2,0.46){\makebox(0,0)[t]{\scriptsize $t_0$}}
	\put(3.2,0.5){\makebox(0,0)[r]{\scriptsize $s_1$}}
	\put(3.2,1){\makebox(0,0)[r]{\scriptsize $s_2$}}
	\put(4.2,-0.01){\makebox(0,0)[t]{\scriptsize $t_1$}}
	\put(4.24,0){\makebox(0,0)[l]{\scriptsize $s_0$}}
	\put(4.24,0.5){\makebox(0,0)[l]{\scriptsize $s_1$}}
	\put(3.18,0.7){\makebox(0,0)[r]{\scriptsize $s$}}
	\put(3.22,0.66){\makebox(0,0)[tl]{\scriptsize $a(s)$}}
	\put(3.8,0.66){\makebox(0,0)[t]{\scriptsize $b(s)$}}
\end{picture}
\end{center}
\end{figure}
		
	To inspect firstly the parallelograms, we consider without loss of generality \linebreak $P_1=P\big((0,s_1),(t_1,0),(t_1,s_1),(0,s_2)\big)$. For each $s\in [0,s_2]$, the line segment $[0,1]\times \{s\}$ intersects $P_1$ in two points $(a(s),s)$ and $(b(s),s)$ with $a(s)\leq b(s)$ that coincide if $s=0$ or $s=s_2$. The hereby defined maps $a,b\colon [0,s_2]\rightarrow [0,1]$ are continuous. The filled parallelogram $\overline P_1$ is contained in some $H^{-1}(W_x)$. We extend $H^\prime|_B$ to $\overline P_1$ by inserting suitable geodesics, i.e.
	$$(H^\prime|_{\overline P_1})(t,s):=\left\{\begin{array}{ll}
			\exp\Big(\frac{t-a(s)}{b(s)-a(s)}\Phi_x^{-1}\big(H^\prime|_B(a(s),s),H^\prime|_B(b(s),s)\big)\Big)
			& \mbox{if $s\neq 0$ and $s\neq s_2$} \\
			H^\prime|_B(t,s) & \mbox{if $(t,s)\in\{(t_1,0),(0,s_2)\}$}.
		\end{array}\right.$$
	It is clear that $H^\prime$ is not redefined on the boundary $P_1$. To verify the continuity, we have to examine more closely the vertices $v\in\{(t_1,0), (0,s_2)\}$: As
	$$\lim_{(t,s)\rightarrow v} \Phi_x^{-1}\big(H^\prime|_B(a(s),s),H^\prime|_B(b(s),s)\big) \ =\ \Phi_x^{-1}\big(H^\prime|_B(v),H^\prime|_B(v)\big) \ =\ 0_{H^\prime|_B(v)}$$
	and as $\frac{t-a(s)}{b(s)-a(s)}\in [0,1]$ is bounded,
	we obtain
	$$\lim_{(t,s)\rightarrow v} H^\prime|_{\overline P_1}(t,s) \ =\ \exp(0_{H^\prime|_B(v)}) \ =\ H^\prime|_B(v).$$
	
	Inspecting the triangles $\Delta_0$ and $\Delta_1$ works quite similarly. An additional subtlety is to check that the geodesics $(H^\prime|_B)(\cdot,0)=\alpha_0|_{[0,t_1]}$ and $(H^\prime|_B)(\cdot,1)=\alpha_1|_{[0,t_1]}$ agree with
	$$t\mapsto \exp\Big(\frac{t}{t_1}\Phi_{x_0}^{-1}\big(H^\prime|_B(0,0),H^\prime|_B(t_1,0)\big)\Big) \quad\mbox{and}\quad t\mapsto \exp\Big(\frac{t}{t_1}\Phi_{x_1}^{-1}\big(H^\prime|_B(0,1),H^\prime|_B(t_1,1)\big)\Big),$$
	respectively, where the filled triangles are contained in some $H^{-1}(W_{x_0})$ and $H^{-1}(W_{x_1})$, respectively. As the geodesics lie in $W_{x_0}$ and $W_{x_1}$, this is ensured by Lemma~\ref{lem:UPhiW}(2).
\end{proof}
\begin{corollary}\label{cor:geodesicHomotopy}
	If $M$ is simply connected, then, given any piecewise geodesics $\alpha_0,\alpha_1\colon$\linebreak $[0,1]\rightarrow M$ with $\alpha_0(0)=\alpha_1(0)$ and $\alpha_0(1)=\alpha_1(1)$, there is a homotopy $H\colon [0,1]\times [0,1]\rightarrow M$ between them such that each curve $\alpha_s:= H(\cdot,s)$ is a piecewise geodesic.
\end{corollary}
\subsection{The Images of Piecewise Geodesics under Maps between Tangent Spaces}
\label{sec:imagesOfPiecewGeod}
Let $(M_1,\nabla\!_1)$ and $(M_2,\nabla\!_2)$ be affine Banach manifolds where the latter is geodesically complete. Given a piecewise geodesic $\alpha$ described by the data
$$(d,\alpha(d);t_0\leq t_1\leq\ldots\leq t_n;v_1,\dots,v_n)$$
and given a continuous linear map $A\colon T_{\alpha(d)}M_1\rightarrow T_yM_2$ (with $y\in M_2$), we define $(d,A)_\ast\alpha$ as the piecewise geodesic given by the data
$$(d,y;t_0\leq t_1\leq\ldots\leq t_n;A(v_1),\dots,A(v_n)).$$
It is easy to see that this is well-defined, i.e., regardless of the partition of the domain of $\alpha$.
We shall mostly write $A_\ast\alpha$ instead of $(d,A)_\ast\alpha$, if the distinguished time $d$ is obvious by the context, e.g., if we use indexed maps like $A_d$.
\begin{definition}\label{def:A_t^alpha}
	Given a piecewise geodesic $\alpha\colon I\rightarrow M_1$ and a continuous linear map\linebreak $A_d\colon T_{\alpha(d)}M_1\rightarrow T_yM_2$, we put $$A_t:=P_{d}^{t}((A_d)_\ast\alpha)\circ A_d\circ P_{t}^{d}(\alpha)\colon T_{\alpha(t)}M_1\rightarrow T_{((A_d)_\ast\alpha)(t)}M_2$$
	for all $t\in I$. To avoid confusion, we sometimes write $A_t^\alpha$ instead of $A_t$.
\end{definition}
\begin{lemma}
	Given a piecewise geodesic $\alpha\colon I\rightarrow M_1$ and a continuous linear map\linebreak $A_d\colon T_{\alpha(d)}M_1\rightarrow T_yM_2$, we have:
	\begin{enumerate}
		\item[\rm (1)]
		$(A_t)_\ast\alpha \ =\ (A_d)_\ast\alpha$ for all $t\in I$.
		\item[\rm (2)]
		$A_t \ =\ P_{s}^{t}((A_s)_\ast\alpha)\circ A_s\circ P_{t}^{s}(\alpha)$ for all $t,s\in I$.
	\end{enumerate}
\end{lemma}
\begin{proof}
	(1) Let $\alpha$ be described by $(d,\alpha(d);t_0\leq t_1\leq\ldots\leq t_n;v_1,\dots,v_n)$. For the sake of readability, we shall omit the time points $t_0,\dots,t_n$.
	Then $(A_d)_\ast\alpha$ can be described by $(d,y;;A_d(v_1),\dots,A_d(v_n))$ or also by
	$$\big(t,((A_d)_\ast\alpha)(t);;P_{d}^{t}((A_d)_\ast\alpha)(A_d(v_1)),\dots,P_{d}^{t}((A_d)_\ast\alpha)(A_d(v_n))\big)$$
	(cf.\ Lemma~\ref{lem:changeOfData}). By change of data, we can describe $\alpha$ by $(t,\alpha(t);;P_{d}^t(\alpha)(v_1),\dots,P_{d}^t(\alpha)(v_n))$, so that $(A_t)_\ast\alpha$ is given by
	$$\big(t,((A_d)_\ast\alpha)(t);;A_t(P_{d}^t(\alpha)(v_1)),\dots,A_t(P_{d}^t(\alpha)(v_n))\big).$$
	By Definition~\ref{def:A_t^alpha}, these data of $(A_d)_\ast\alpha$ and $(A_t)_\ast\alpha$ agree.
	
	(2) We have
	\begin{eqnarray*}
		P_{s}^{t}((A_s)_\ast\alpha)\circ A_s\circ P_{t}^{s}(\alpha)
		&=& P_{s}^{t}((A_s)_\ast\alpha)\circ (P_{d}^{s}((A_d)_\ast\alpha)\circ A_d\circ P_{s}^{d}(\alpha)) \circ P_{t}^{s}(\alpha) \\
		&\stackrel{\mbox{\scriptsize (1)}}{=}&
		P_{d}^{t}((A_d)_\ast\alpha)\circ A_d\circ P_{t}^{d}(\alpha) \ =\ A_t.
	\end{eqnarray*}
\end{proof}
\begin{corollary}\label{cor:sameCollection}
	If we start in the situation of Definition~\ref{def:A_t^alpha} with the map $A_{d^\prime}$ (induced by $A_d$) for another distinguished time $d^\prime$, then we get the same collection of maps $(A_t)_{t\in I}$.
\end{corollary}
\begin{lemma}\label{lem:collectionT_alpha(t)f}
	Given a piecewise geodesic $\alpha\colon I\rightarrow M_1$ and an affine map $f\colon M_1\rightarrow M_2$, we have
	$f\circ\alpha = (T_{\alpha(d)}f)_\ast\alpha$
	for all $d\in I$. The collection $(T_{\alpha(t)}f)_{t\in I}$ is a collection in the sense of Corollary~\ref{cor:sameCollection}.
\end{lemma}
\begin{proof}
	Let $\alpha=(\alpha_1,\dots,\alpha_n)$ be described by $(d,\alpha(d);t_0\leq t_1\leq\ldots\leq t_n;v_1,\dots,v_n)$. Then $(T_{\alpha(d)}f)_\ast\alpha$ is given by the data
	$$(d,f(\alpha(d));t_0\leq t_1\leq\ldots\leq t_n;(T_{\alpha(d)}f)(v_1),\dots,(T_{\alpha(d)}f)(v_n)).$$
	To see that this data also describes $f\circ\alpha$, we have to check
	$$(f\circ\alpha_i)^\prime(t)=P_{d}^{t}(f\circ\alpha)((T_{\alpha(d)}f)(v_i))\quad\mbox{for all $t \in I_i=[t_{i-1},t_i]$ with $i=1,\dots,n$}.$$
	The latter is true, since we have
	$$(f\circ\alpha_i)^\prime(t) \ =\ Tf(\alpha_i^\prime(t)) \ =\ Tf(P_{d}^t(\alpha)(v_i)) \ =\ P_{d}^t(f\circ\alpha)(Tf(v_i)),$$
	affine maps being compatible with parallel transport. Hence, we have $f\circ\alpha = (T_{\alpha(d)}f)_\ast\alpha$.
	
	As we have $$T_{\alpha(t)}f \ =\ P_{s}^t(f\circ\alpha)\circ T_{\alpha(s)}f \circ P_{t}^s(\alpha) \ =\ P_{s}^t((T_{\alpha(s)}f)_\ast\alpha)\circ T_{\alpha(s)}f \circ P_{t}^s(\alpha)$$
	for all $t,s\in I$, the collection $(T_{\alpha(t)}f)_{t\in I}$ is a collection in the sense of Corollary~\ref{cor:sameCollection}.
\end{proof}
\subsection{Integration of Maps between Tangent Spaces}
\label{sec:subsectionIntegrationToAffMaps}
\begin{proposition} \label{prop:integrabilityToAffineMaps}
	Let $(M_1,\nabla\!_1,b_1)$ and $(M_2,\nabla\!_2,b_2)$ be affine Banach manifolds with base points where the former is connected and the latter is geodesically complete. Given a continuous linear map $A\colon T_{b_1}M_1\rightarrow T_{b_2}M_2$, the following are equivalent:
	\begin{enumerate}
		\item[\rm (a)]
		Given any piecewise geodesics $\alpha,\beta\colon [0,1]\rightarrow M_1$ with $\alpha(0)=\beta(0)=b_1$ and $\alpha(1)=\beta(1)$,
		we have $A^\alpha_1=A^\beta_1$ for $A^\alpha_0:=A^\beta_0:=A$. Further, each map $A^\alpha_1$ is locally integrable to an affine map $g\colon U_{\alpha(1)}\rightarrow M_2$ with an open neighborhood $U_{\alpha(1)}$ of $\alpha(1)$ in the sense that $T_{\alpha(1)}g=A^\alpha_1$.
		\item[\rm (b)]
		There exists an affine map $f\colon M_1\rightarrow M_2$ that satisfies $T_{b_1}f=A$.
	\end{enumerate}
\end{proposition}
\begin{proof}
	(b)$\Rightarrow$(a): By Lemma~\ref{lem:collectionT_alpha(t)f}, we have $A^\alpha_1=T_{\alpha(1)}f$ and $A^\beta_1=T_{\beta(1)}f$, so that $A^\alpha_1=A^\beta_1$.
	
	(a)$\Rightarrow$(b): We define the map $f\colon M_1\rightarrow  M_2$ by
	\begin{equation} \label{eqn:f(alpha(1))}
		f(\alpha(1)):=((A^\alpha_0)_\ast\alpha)(1) \quad\mbox{for all $\alpha\colon [0,1]\rightarrow M_1$ with $\alpha(0)=b_1$ and $A^\alpha_0:=A$.}
	\end{equation}
	Due to Lemma~\ref{lem:geodesicConnection}, $f$ is defined on all of $M$. It is well-defined, as	$$\alpha(1)=\beta(1) \ \stackrel{\mbox{\scriptsize (a)}}{\Rightarrow} \ A^\alpha_1=A^\beta_1 \ \Rightarrow \ ((A^\alpha_0)_\ast\alpha)(1) = ((A^\beta_0)_\ast\beta)(1)$$
	for all such piecewise geodesics $\alpha$ and $\beta$.
	
	For each $x\in M$, we define the map
	$A_x\colon T_xM_1\rightarrow T_{f(x)}M_2$ by $A_x:=A^\alpha_1$, where\linebreak $\alpha\colon [0,1]\rightarrow M_1$ is any piecewise geodesic that joins $b_1$ with $x$ and where $A^\alpha_0:=A$. It is well-defined by (a). It is clear that $A_{b_1}=A$.
	
	We claim that
	\begin{equation*}
		f(\alpha(1))=((A^\alpha_0)_\ast\alpha)(1) \quad\mbox{for all $\alpha\colon [0,1]\rightarrow M_1$ with $\alpha(0)=x$ and $A^\alpha_0:=A_x$.}
	\end{equation*}
	To see this, let $\beta\colon [-1,0]\rightarrow M_1$ be a piecewise geodesic that joins $b_1$ with $x$ and let $(\beta,\alpha)\colon [-1,1] \rightarrow M_1$ be the piecewise geodesic that is obtained by composing $\beta$ and $\alpha$. Putting $A^{(\beta,\alpha)}_{-1}:=A$ and applying the known formulas after reparametrization, we obtain
	\begin{eqnarray*}
		f(\alpha(1)) &=& f((\beta,\alpha)(1)) \ \stackrel{\mbox{\scriptsize (\ref{eqn:f(alpha(1))})}}{=}\ ((A^{(\beta,\alpha)}_{-1})_\ast(\beta,\alpha))(1) \ =\ ((A^{(\beta,\alpha)}_{0})_\ast(\beta,\alpha))(1)\\
		&\stackrel{!}{=}& ((A^\alpha_{0})_\ast(\beta,\alpha))(1) \ =\ ((A^\alpha_{0})_\ast\alpha))(1),
	\end{eqnarray*}
	since
	$A^{(\beta,\alpha)}_0 = A_x = A^\alpha_0$.
		
	We shall show that $f$ is smooth and affine by working locally: Given any $x\in M$, we integrate $A_x$ to a smooth affine map $g\colon U_x \rightarrow M_2$ with $T_xg=A_x$, where $U_x$ is some connected neighborhood of $x$. We claim that $g$ is a restriction of $f$. Indeed, given any $y\in U_x$, we consider some piecewise geodesic $\alpha\colon [0,1]\rightarrow U_x$ that joins $x$ with $y$, so that
	$g(y)=(g\circ\alpha)(1)$ and $f(y)=f(\alpha(1))=((A^\alpha_0)_\ast\alpha)(1)$ with $A^\alpha_0:=A_x=T_xg$. By Lemma~\ref{lem:collectionT_alpha(t)f}, we then have $g=f|_{U_x}$.
	
	In particular, we have $T_{b_1}f=A_{b_1}=A$.
\end{proof}
\begin{definition}\label{def:uniformlyLocallyIntegrable}
	Let $(M_1,\nabla\!_1)$ and $(M_2,\nabla\!_2)$ be affine Banach manifolds where the latter is geodesically complete and let $\alpha\colon I \rightarrow M_1$ be a piecewise geodesic. Given a collection $(A_t)_{t\in I}$ of continuous linear maps in the sense of Corollary~\ref{cor:sameCollection}, we call it \emph{uniformly locally integrable} if there are open neighborhoods $U_t$ of $\alpha(t)$ and affine maps $f_t\colon U_t\rightarrow M_2$ with $T_{\alpha(t)}f_t=A_t$ for all $t\in I$ and some decomposition $\alpha=(\alpha_1,\dots,\alpha_n)$ with geodesics $\alpha_i\colon I_i\rightarrow M_1$ (cf.\ Section~\ref{sec:piecewiseGeodesics}) such that
	$U_t \supseteq \im(\alpha_i)$ for all $i=1,\dots, n$ and $t\in I_i$.
\end{definition}
\begin{lemma} \label{lem:f_s=f_t}
	Let $(M_1,\nabla\!_1)$ and $(M_2,\nabla\!_2)$ be affine Banach manifolds where the latter is geodesically complete and let $\alpha\colon I \rightarrow M_1$ be a piecewise geodesic. Given a uniformly locally integrable collection $(A_t)_{t\in I}$, we have (with respect to the denotations of Definition~\ref{def:uniformlyLocallyIntegrable})
	$$U_s\cap U_t\supseteq \im(\alpha_i) \quad\mbox{and}\quad f_s|_{(U_s\cap U_t)_c}= f_t|_{(U_s\cap U_t)_c}$$
	for all $i=1,\dots, n$ and $[s,t]\subseteq I_i$, where $(U_s\cap U_t)_c$ is the (open) connected component of $(U_s\cap U_t)$ that contains $\im(\alpha_i)$.
\end{lemma}
\begin{proof}
	By Section~\ref{sec:affineMaps}, it suffices to check $T_{\alpha_i(s)}f_s=T_{\alpha_i(s)}f_t$. The collections $(T_{\alpha_i(d)}f_t)_{d\in I_i}$ (cf.\ Lemma~\ref{lem:collectionT_alpha(t)f}) and $(A_d)_{d\in I_i}$ agree, as in particular $T_{\alpha_i(t)}f_t=A_t$. Hence, we have $T_{\alpha_i(s)}f_t=A_s=T_{\alpha_i(s)}f_s$.
\end{proof}
\begin{lemma}
	Let $(M_1,\nabla\!_1)$ and $(M_2,\nabla\!_2)$ be affine Banach manifolds where the latter is geodesically complete and let $\alpha\colon I \rightarrow M_1$ be a piecewise geodesic. Given a collection $(A_t)_{t\in I}$ of continuous linear maps in the sense of Corollary~\ref{cor:sameCollection}, it is uniformly locally integrable if and only if each map $A_t$ is locally integrable.
\end{lemma}
\begin{proof}
	Assume that there are open neighborhoods $V_t$ of $\alpha(t)$ and affine maps $g_t\colon V_t\rightarrow M_2$ with $T_{\alpha(t)}g_t=A_t$ for all $t\in I$. They form a covering $(V_t)_{t\in I}$ of the image $\im(\alpha)$ of $\alpha$. By a compactness argument, there is a finite subcover $(V_t)_{t\in F}$, some decomposition $\alpha=(\alpha_1,\dots,\alpha_n)$ with geodesics $\alpha_i\colon I_i\rightarrow M_1$ (cf.\ Section~\ref{sec:piecewiseGeodesics}) and a map $j\colon I\rightarrow F$ such that for each $t\in I_i$, $i=1,\ldots,n$, the convex hull $C_{i,t}:=\conv(I_i\cup \{j(t)\})$ satisfies $V_{j(t)}\supseteq \alpha(C_{i,t})$.\footnote{
	Let $(W_t)_{t\in I}$ be a covering of $I$ where each $W_t$ is an open convex neighborhood of $t$ such that $\alpha(W_t)\subseteq V_t$. We choose a finite subset $F$ of $I$ such that for all adjacent points $s$ and $t$ in $F$, the neighborhoods $W_t$ and $W_{s}$ intersect. We decompose $I$ by choosing two points in each such intersection. Then we can define a map $j\colon I \rightarrow F$ such that for each $t\in I$, the set $W_{j(t)}$ covers the subinterval of $I$ containing $t$. Note in this situation that the dividing points belong to two intervals. Establishing a refinement ensures that we decompose $\alpha$ into geodesic pieces.}
	For each $t\in I$, we then define $U_t:=V_{j(t)}$ and $f_t:=g_{j(t)}$ and claim that $T_{\alpha(t)}f_t=A_t$. Indeed, the collections $(T_{\alpha(d)}f_t)_{d\in C_{i,t}}$ (cf.\ Lemma~\ref{lem:collectionT_alpha(t)f}) and $(A_d)_{d\in C_{i,t}}$ agree, as in particular $T_{\alpha(j(t))}f_t=T_{\alpha(j(t))}g_{j(t)}=A_{j(t)}$.
\end{proof}
\begin{proposition} \label{prop:independantA_1}
	Let $(M_1,\nabla\!_1,b_1)$ and $(M_2,\nabla\!_2,b_2)$ be affine Banach manifolds with base points where the latter is geodesically complete. Given a continuous linear map \linebreak$A\colon T_{b_1}M_1\rightarrow T_{b_2}M_2$, two piecewise geodesics $\alpha_0,\alpha_1\colon [0,1]\rightarrow M_1$ with $\alpha_0(0)=b_1=\alpha_1(0)$ and $\alpha_0(1)=\alpha_1(1)$ and a homotopy $H\colon [0,1]\times [0,1]\rightarrow M_1$ between $\alpha_0$ and $\alpha_1$ such that each curve $\alpha_s:=H(\cdot,s)$ is a piecewise geodesic, we put $A^{\alpha_s}_0:=A$ for each $s\in [0,1]$.
	If for each $s\in [0,1]$ and $t\in I$, the map $A^{\alpha_s}_t$ is locally integrable, then $A^{\alpha_s}_1$ is independent of $s$.
\end{proposition}
\begin{proof}
	We shall prove the assertion in three steps.
	
	{\bf Step 1:} \emph{For fixed $s\in [0,1]$, we consider some decomposition $\alpha_s=(\alpha_{s,1},\dots,\alpha_{s,n})$ with geodesics $\alpha_{s,i}\colon I_i=[t_{i-1},t_i]\rightarrow M_1$
	and affine maps $f_t\colon U_t\rightarrow M_2$ as in Definition~\ref{def:uniformlyLocallyIntegrable}. Then there is an open neighborhood $V^s$ of $s$ in $[0,1]$ such that $H([t_{i-1},t_i]\times V^s)\subseteq (U_{t_{i-1}}\cap U_{t_i})_c$ for all $i=1,\dots,n$.}
	
	To see this, we put $t_{-1}:=t_0=0$ and $t_{n+1}:=t_n=1$ for the sake of readability and observe that for each $i=0,\dots,n$, the open preimage $H^{-1}(U_{t_i})$ contains $[t_{i-1},t_{i+1}]\times \{s\}$, so that, by a compactness argument, there is an open convex neighborhood $V^s_i$ of $s$ in $[0,1]$ such that $H([t_{i-1},t_{i+1}]\times V^s_i)\subseteq U_{t_i}$. By putting $V^s:=\cap_{i=0}^n V^s_i$, we obtain $H([t_{i-1},t_i]\times V^s)\subseteq U_{t_{i-1}}\cap U_{t_i}$ and hence $H([t_{i-1},t_i]\times V^s)\subseteq (U_{t_{i-1}}\cap U_{t_i})_c$ for all $i=1,\dots, n$.
	
	{\bf Step 2:} \emph{For fixed $s\in [0,1]$, we have $A^{\alpha_{s^\prime}}_1=A^{\alpha_s}_1$ for all $s^\prime\in V^s$.}
	To see this, we note $A^{\alpha_s}_1=T_{\alpha_{s}(1)}f_1=T_{\alpha_{s^\prime}(1)}f_1$ and prove the more general assertion $A^{\alpha_{s^\prime}}_{t_i}=T_{\alpha_{s^\prime}(t_i)}f_{t_i}$ for all $i=0,\dots,n$. For $i=0$ this is trivial, as $A^{\alpha_{s^\prime}}_0=A=T_{b_1}f_0$.
	It suffices to verify the step $i\rightarrow i+1$ in the sense that we assume the assertion for given $i<n$ and deduce it for $i+1$. We first note that
	$\alpha_{s^\prime}([t_i,t_{i+1}])=H([t_i,t_{i+1}]\times \{s^\prime\})\subseteq (U_{t_i}\cap U_{t_{i+1}})_c\subseteq U_{t_i}$.
	The collections $(T_{\alpha_{s^\prime}(d)}f_{t_i})_{d\in [t_i,t_{i+1}]}$ (cf.\ Lemma~\ref{lem:collectionT_alpha(t)f}) and $(A^{\alpha_{s^\prime}}_d)_{d\in [t_i,t_{i+1}]}$ agree, as in particular $T_{\alpha_{s^\prime}(t_i)}f_{t_i}=A^{\alpha_{s^\prime}}_{t_i}$, by assumption. Hence, we have $T_{\alpha_{s^\prime}(t_{i+1})}f_{t_i}=A^{\alpha_{s^\prime}}_{t_{i+1}}$. The maps $f_{t_i}$ and $f_{t_{i+1}}$ agreeing on $(U_{t_i}\cap U_{t_{i+1}})_c$ (cf.\ Lemma~\ref{lem:f_s=f_t}), we obtain $T_{\alpha_{s^\prime}(t_{i+1})}f_{t_{i+1}}=A^{\alpha_{s^\prime}}_{t_{i+1}}$.
	
	{\bf Step 3:} \emph{Observe that $A^{\alpha_s}_1$ is independent of $s$.}
	To see this, let $S$ be the set of all $s\in [0,1]$ for which $A^{\alpha_s}_1=A^{\alpha_0}_1$. We have to show that $S$ is all of $[0,1]$. Due to Step~2, the non-empty set $S$ and its complement $S^c$ both are open in $[0,1]$, so that $S=[0,1]$ by the connectedness of $[0,1]$.
\end{proof}
\begin{theorem}
\label{th:cartan-ambrose-hicks}
	Let $(M_1,\nabla\!_1,b_1)$ and $(M_2,\nabla\!_2,b_2)$ be affine Banach manifolds with base points where the former is 1-connected, i.e., connected and simply connected, and the latter is geodesically complete. Given a continuous linear map $A\colon T_{b_1}M_1\rightarrow T_{b_2}M_2$, the following are equivalent:
	\begin{enumerate}
		\item[\rm (a)]
		Given any piecewise geodesic $\alpha\colon [0,1] \rightarrow M_1$ with $\alpha(0)=b_1$, the map $A_1$ with $A_0:=A$ (cf.\ Definition~\ref{def:A_t^alpha}) is locally integrable.
		\item[\rm (b)] There exists an affine map $f\colon M_1\rightarrow M_2$ that satisfies $T_{b_1}f=A$.
	\end{enumerate}
\end{theorem}
\begin{proof}
	(b)$\Rightarrow$(a) is clear (cf.\ Lemma~\ref{lem:collectionT_alpha(t)f}) and (a)$\Rightarrow$(b) follows immediately by Proposition~\ref{prop:integrabilityToAffineMaps}, Proposition~\ref{prop:independantA_1} and Corollary~\ref{cor:geodesicHomotopy}.
\end{proof}
\begin{lemma}\label{lem:collectionPreservesTorsionAndCurvature}
	Let $(M_1,\nabla\!_1)$ and $(M_2,\nabla\!_2)$ be affine Banach manifolds where the latter is geodesically complete and where the torsion and curvature tensors are parallel.
	Given a piecewise geodesic $\alpha\colon I\rightarrow M_1$ and a continuous linear map $A_d\colon T_{\alpha(d)}M_1\rightarrow T_yM_2$ (with $d\in I$ and $y\in M_2$) that intertwines
	the torsion and curvature tensors in $x:=\alpha(d)$ and $y$, i.e.,
	$$A_d\circ (\Tor_1)_{x} = (\Tor_2)_{y}\circ A_d^2 \quad\mbox{and}\quad A_d\circ (R_1)_{x} = (R_2)_{y}\circ A_d^3,$$
	then each map of the collection $(A_t)_{t\in I}$ intertwines the torsion and curvature tensors in $\alpha(t)$ and $((A_d)_\ast\alpha)(t)$.
\end{lemma}
\begin{proof}
	For each $t\in I$, we have $A_t=P_{d}^{t}((A_d)_\ast\alpha)\circ A_d\circ P_{t}^{d}(\alpha)$. Therefore, it suffices to show that parallel transport along curves preserves torsion and curvature, but this is true, as these tensors are assumed to be parallel (cf.\ Section~\ref{sec:affineConnections}).
\end{proof}
\begin{theorem}[{Integrability Theorem}]\label{th:globalIntegrability}
	Let $(M_1,\nabla\!_1,b_1)$ and $(M_2,\nabla\!_2,b_2)$ be affine Banach manifolds with base points where the former is 1-connected  and the latter is geodesically complete and where the torsion and curvature tensors are parallel. For every continuous linear map $A\colon T_{b_1}M_1\rightarrow T_{b_2}M_2$ that intertwines the torsion and curvature tensors in the base points in the sense that
	$$A\circ (\Tor_1)_{b_1} = (\Tor_2)_{b_2}\circ A^2 \quad\mbox{and}\quad A \circ (R_1)_{b_1} = (R_2)_{b_2}\circ A^3,$$
	there exists a unique affine map $f\colon M_1\rightarrow M_2$ that satisfies $T_{b_1}f=A$.
\end{theorem}
\begin{proof}
	The theorem follows by Theorem~\ref{th:cartan-ambrose-hicks}, Lemma~\ref{lem:collectionPreservesTorsionAndCurvature} and Corollary~\ref{cor:extensionToAffMaps}. Cf.\ Section~\ref{sec:affineMaps} for the uniqueness assertion.
\end{proof}
\begin{remark}
	This theorem generalizes a special case of the theorem of Cartan--Ambrose--Hicks to the Banach case, where furthermore the map $A$ is not supposed to be an isomorphism.
	In \cite{Amb56}, W.~Ambrose gives a theorem about the integration of an isometric isomorphism between tangent spaces of complete 1-connected Riemannian manifolds in the finite-dimensional case. As a sequel to it, in \cite{Hic59}, N.~Hicks deals with the case of affine manifolds.
	Their work is based on Cartan's work on frames and connection 1-forms (cf.\ \cite{Car46}).
\end{remark}
%
%
%
%
%
%
\section{Banach Symmetric Spaces} \label{sec:symmetricSpaces}
In \cite{Loo69}, O.~Loos defines symmetric spaces by means of a multiplication map $\mu$ on a finite-dimensional manifold $M$,
where each left multiplication map $\mu_x$ (with $x\in M$) is an involutive automorphism of $(M,\mu)$ with the isolated fixed point $x$.
In this section, we shall deal with symmetric spaces in the context of Banach manifolds.

For each point of a symmetric space, there is a natural involutive automorphism of the Lie algebra $\Der(M,\mu)$ of derivations that is induced by the symmetry at this point. This provides an additional structure on the tangent space, namely a Lie triple system. Furthermore, there is a functor $\Lts$ from the category of pointed symmetric spaces to the category of Lie triple systems.

The main results of this section are an integrability theorem about morphisms of Lie triple systems and the fact that the automorphism group of a connected symmetric space $M$ can be turned into a Banach--Lie group acting transitively on $M$. As a consequence, we shall see that connected symmetric spaces are homogeneous.

To obtain these results, we shall equip a symmetric space with a canonical affine connection encoding the symmetric space structure in the sense that it has the same automorphisms. Observing that symmetric spaces are torsionfree and have parallel curvature and that morphisms of Lie triple systems of pointed symmetric spaces are just the curvature preserving maps, we can apply the preceding results.

\subsection{Banach Symmetric Spaces}
\label{sec:subsectionSymmetricSpaces}
	A \emph{Banach symmetric space}, simply called a \emph{symmetric space}, is a smooth Banach manifold $M$ with a smooth multiplication $\mu \colon M\times M \rightarrow M$, written as $\mu(x,y)=x \cdot y$, such that for all $x,y,z\in M$, writing $\mu_x:=\mu(x,\cdot)$ for the left multiplication, the following properties hold:
	\begin{enumerate}
		\item[(S1)] $x \cdot x =x$, i.e., $\mu_x(x)=x$.
		\item[(S2)] $x\cdot(x\cdot y) = y$, i.e., $\mu_x^2 = \id_M$.
		\item[(S3)] $x\cdot(y \cdot z) = (x \cdot y) \cdot (x \cdot z)$,
			i.e., $\mu_x(y\cdot z)=\mu_x(y)\cdot\mu_x(z)$.
		\item[(S4)] Every $x$ has a neighborhood $U$ such that $x\cdot y=y$ implies $y=x$ for all $y\in U$, i.e., $x$ is an isolated fixed point of $\mu_x$.
	\end{enumerate}
	We mention (but will not make use of this) that (S4) can be replaced by the condition
	\begin{enumerate}
		\item[(S4')] $T_x\mu_x=-\id_{T_xM}$
	\end{enumerate}
	(cf.\ \cite[Lem.~3.2]{Nee02}).\footnote{The proof of the cited lemma is incorrect, but reparable: Considering a local representation $\mu_x^\varphi\colon V\rightarrow V$ of $\mu_x$ (for a suitable chart $\varphi\colon U\rightarrow V\subseteq E$), we can turn $V$ into an affine manifold making $\mu_x^\varphi$ an affine involutive diffeomorphism. Since the set of all affine connections on $V$ is an affine space (cf.\ \cite[10.4]{Ber08}), such an affine connection can be obtained as the midpoint of the trivial flat connection on $V$ and its pushforward along $\mu_x^\varphi$. Then we have $\mu_x^\varphi\circ\exp_{\varphi(x)} = \exp_{\varphi(x)}\circ T_{\varphi(x)}\mu_x^\varphi$, entailing that $\varphi(x)$ is an isolated fixed point of $\mu_x^\varphi$ if and only if $0\in T_{\varphi(x)}V$ is an isolated fixed point of $T_{\varphi(x)}\mu_x^\varphi$, which is equivalent to $T_{\varphi(x)}\mu_x^\varphi=-\id_{T_{\varphi(x)}\mu_x^\varphi}$.}
	
	A \emph{morphism between symmetric spaces $(M_1,\mu_1)$ and $(M_2,\mu_2)$} is a smooth map\linebreak $f\colon M_1 \rightarrow M_2$ such that
	$$f\circ\mu_1 \ =\ \mu_2\circ(f\times f),$$
	i.e., $f(x\cdot y)= f(x)\cdot f(y)$ for all $x,y\in M_1$.
	The class of symmetric spaces and the class of morphisms between them form a category, so that isomorphisms and automorphisms are defined as usual.
	For each $x\in M$, the map $\mu_x$ is called the \emph{symmetry around x}. Obviously, it is an involutive automorphism of $(M,\mu)$.
	A \emph{pointed symmetric space} is a triple $(M,\mu,b)$ consisting of a symmetric space $(M,\mu)$ and a point $b\in M$ called the \emph{base point}. A morphism $f$ between pointed symmetric spaces $(M_1,\mu_1,b_1)$ and $(M_2,\mu_2,b_2)$ is required to be base-point preserving, i.e., $f(b_1)=b_2$. We call it a \emph{morphism of pointed symmetric spaces}. Note that $\mu_b\in\Aut(M,\mu,b)$.
%
%
\begin{example}[Lie groups] \label{ex:lieGroups}
	Let $G$ be a Banach--Lie group. Then the manifold $G$ together with the map $\mu(g,h):=gh^{-1}g$, where $gh$ denotes the product in $G$, is a symmetric space.
	In particular, if $G$ is a Banach space, it becomes a symmetric space with the multiplication $\mu(g,h):=2g-h$. For further details, see \cite[pp.~65-66]{Loo69}, which carries over to the Banach case.
\end{example}
\begin{example}[Homogeneous spaces] \label{ex:homogeneousSpaces}
	Given a Banach--Lie group $G$ and an involutive automorphism $\sigma$ of $G$, let $G^\sigma:=\{g\in G\colon \sigma(g)=g\}$ be the subgroup of $\sigma$-fixed points and $K\subseteq G^\sigma$ an open subgroup. The quotient space $M:=G/K$ carries the structure of a Banach manifold such that the quotient map $q\colon G\rightarrow G/K$ is a submersion. It can be equipped with a multiplication
	$$\mu(gK,hK):=g\sigma(g)^{-1}\sigma(h)K$$
	that turns it into a symmetric space. For further details, cf.\ \cite[Ex.~3.9]{Nee02}.
\end{example}
\begin{example}[Quadrics]
	Let $E$ be a Banach space with a continuous symmetric bilinear form
	$\langle\cdot,\cdot\rangle\colon E\times E \rightarrow \RR$. Given $a\in\RR^{\times}$, the quadric $S:=\{x\in E\colon \langle x,x\rangle = a\}$ is a (split) submanifold of $E$. Indeed, for any $x\in E$, the derivative $df(x)$ of the smooth map \linebreak $f\colon x\mapsto \langle x,x\rangle$ is given by
	$$df(x)=2\langle x,\cdot\rangle$$
	and is therefore surjective if $x\in S$. Furthermore, its kernel splits for $x\in S$, since $E = x^\bot \oplus \RR x$ according to the decomposition
	$$E \ni y = \Big(y-\frac{\langle x,y\rangle}{\langle x,x\rangle}x\Big)
				+ \frac{\langle x,y\rangle}{\langle x,x\rangle}x,$$
	where $x^\bot$ denotes the kernel of $\langle x,\cdot\rangle$.
	Hence the maps $df(x)$ are submersions, so that due to $S=f^{-1}(a)$, the quadric is a submanifold. Note that $x^\bot=\ker(df(x))$ can be identified with the tangent space $T_xS$. For further details, see \cite[II, \S 2]{Lan01} which deals with the Hilbert case.
	
	By equipping $S$ with the multiplication
	$$\mu(x,y):= 2\frac{\langle x,y\rangle}{\langle x,x\rangle}x - y,$$
	we turn it into a symmetric space. The finite-dimensional case can be found in \cite[p.~66]{Loo69}. It carries over to the Banach case without difficulty.
\end{example}
\begin{example}[Spaces of symmetric elements]
\label{ex:symmetricElements}
	Given a Banach--Lie group $G$ and an involutive automorphism $\sigma$ of $G$, the set
	$$G_\sigma:=\{g\in G\colon \sigma(g)=g^{-1}\}$$
	of symmetric elements is a (split) submanifold of $G$. Indeed, let
	$$\varphi:=(\exp|_V^U)^{-1}\colon U \rightarrow V\subseteq L(G)$$ be an exponential chart with $V=-V$. Then for any $g\in G_\sigma$, the chart
	$$\begin{array}{cccc}
		\varphi_g:=\varphi\circ\lambda_{g^{-1}}|_{gU}\colon & gU & \rightarrow & V \\
		& g\exp(x) & \mapsto & x
	\end{array}$$
	maps the intersection $gU\cap G_\sigma$ onto $\{x\in V\colon g\exp(x)\in G_\sigma\}$, 
	which equals
	$$\{x\in V\colon \exp\big((L(\sigma)\circ \Ad(g))(x)\big)=\exp(-x)\},$$
	since
	\begin{eqnarray*}
		\sigma(g\exp(x))=(g\exp(x))^{-1} & \Leftrightarrow & \sigma(g\exp(x)g^{-1})=(\exp(x))^{-1} \\
			&\Leftrightarrow & \exp\big((L(\sigma)\circ\Ad(g))(x)\big) = \exp(-x).
	\end{eqnarray*}
	We make the chart $\varphi_g$ sufficiently small such that $(L(\sigma)\circ\Ad(g))(V)$ lies in the original $V$. Then $\exp\big((L(\sigma)\circ\Ad(g))(x)\big)=\exp(-x)$ is equivalent to $(L(\sigma)\circ\Ad(g))(x)=-x$, so that
	$$\varphi_g(gU\cap G_\sigma)=V\cap\{x\in L(G)\colon (L(\sigma)\circ\Ad(g))(x)=-x\}.$$
	The map $L(\sigma)\circ\Ad(g) = L(\sigma\circ c_g)$ (where $c_g$ denotes the conjugation map) is an involution on $L(G)$ because $\sigma\circ c_g$ is an involution on $G$. Hence, the set $\{x\in L(G)\colon (L(\sigma)\circ\Ad(g))(x)=-x\}$ is a closed subspace of $L(G)$ that is complemented by $\{x\in L(G)\colon (L(\sigma)\circ\Ad(g))(x)=x\}$.
	
	If we view the Lie group $G$ as a symmetric space (cf.\ Example~\ref{ex:lieGroups}), the submanifold $G_\sigma$ inherits this structure, since it is stable under products:
	\begin{eqnarray*}
		g,h\in G_\sigma & \Rightarrow &
			\sigma(gh^{-1}g)= \sigma(g)\sigma(h)^{-1}\sigma(g)=g^{-1}hg^{-1}=(gh^{-1}g)^{-1} \\
		& \Rightarrow & g\cdot h = gh^{-1}g\in G_\sigma
	\end{eqnarray*}
	Hence, $(G_\sigma,\mu)$ with $\mu(g,h):=gh^{-1}g=g\sigma(h)g$ is a symmetric space.
\end{example}

\begin{example}[Spaces of involutions]

	(a) Given a Banach--Lie group $G$, the set
	$$\Invol(G):=\{g\in G\colon g^2=\1\}=\{g\in G\colon g=g^{-1}\}$$
	of involutions is a (split) submanifold of $G$ that inherits from $G$ the structure of a symmetric space with the multiplication map $\mu(g,h):=gh^{-1}g=ghg$ (cf.\ Example~\ref{ex:lieGroups}). This is a special case of Example~\ref{ex:symmetricElements} with $\sigma=\id_G$.

	(b) Given a Banach algebra $A$, the set
	$$\Invol(A):=\{x\in A\colon x^2=\1\}$$
	of involutions coincides with the set of involutions of the open unit group $A^\times$. Hence it carries the structure of a symmetric space (cf.\ (a)). The multiplication $\mu$ is given by $\mu(x,y):= xyx$.
\end{example}
\begin{example}[Grassmannians]
	(a) Given a Banach space $E$, let $\Gr(E)$ denote the set of all closed subspaces that split in $E$. For every subspace $F_2$ of $E$, let $U_{F_2}$ denote the set of all closed subspaces $F_1$ of $E$ that complement $F_2$, i.e., $E=F_1\oplus F_2$. For each such pair $(F_1,F_2)$, we define a bijection
	$$\varphi_{(F_1,F_2)}\colon U_{F_2} \rightarrow \calL(F_1,F_2)$$
	that maps every $F\in U_{F_2}$ to a continuous linear map $f\in \calL(F_1,F_2)$ such that the subspace $F$ of $E=F_1\times F_2$ is the graph of $f$. The charts $\varphi_{(F_1,F_2)}$ turn $\Gr(E)$ into a smooth Banach manifold (cf.\ \cite[5.2.6]{Bou07}).
	
	For each decomposition $E=F_1\oplus F_2$, let $\sigma_{(F_1,F_2)}$ be the reflection in $F_1$, i.e.,
	$$\sigma_{(F_1,F_2)}(x_1\oplus x_2)= x_1 \oplus (-x_2)$$
	for all $x_1\oplus x_2\in F_1\oplus F_2$. If we equip the open submanifold
	$${\cal D} = \{(F_1,F_2)\in \Gr(E)^2\colon E=F_1\oplus F_2\}$$
	of $\Gr(E)^2=\Gr(E)\times\Gr(E)$ with the multiplication
	$$\mu((F_1,F_2),(F_1^\prime,F_2^\prime))
		=(\sigma_{(F_1,F_2)}(F_1^\prime),\sigma_{(F_1,F_2)}(F_2^\prime)),$$
	then it becomes a symmetric space.
	
	(b) Given a scalar product on $E$ that turns it into a Hilbert space, the Grassmannian $\Gr(E)$ can be considered as a submanifold of $\cal D$ according to the embedding $\iota\colon F\mapsto (F,F^\bot)$. Indeed, for each $(F,F^\bot)$, the chart
	$$\varphi_{(F,F^\bot)}\times\varphi_{(F^\bot,F)}\colon
		U_{F^\bot}\times U_F \rightarrow \calL(F,F^\bot)\times \calL(F^\bot,F)$$
	of $\Gr(E)^2$ at $(F,F^\bot)$ maps the set $(U_{F^\bot}\times U_F) \cap \iota (\Gr(E))$ to the closed subspace
	$$\{(f,g)\in\calL(F,F^\bot)\times \calL(F^\bot,F)\colon f+g^\ast=0\},$$
	which is complemented by
	$$\{(f,g)\in\calL(F,F^\bot)\times \calL(F^\bot,F)\colon f-g^\ast=0\}.$$
	Being stable under products, $\Gr(E)$ becomes a symmetric space.
		
	This can be seen also more directly by defining $\sigma_F$ to be the reflection in $F$, i.e.
	$$\sigma_F(x_1\oplus x_2)= x_1 \oplus (-x_2)$$
	for all $x_1\oplus x_2\in F \oplus F^\bot=E$,
	and by defining the multiplication
	$$\mu(F_1,F_2):=\sigma_{F_1}(F_2)$$
	on $\Gr(E)$.
	The finite-dimensional case is dealt with in \cite[pp.~66, 67]{Loo69}.
\end{example}
\begin{remark}
	Given a Banach space $E$, we consider the Banach algebra $\calL(E)$ of endomorphisms of $E$.
	Each involution $A\in\calL(E)$ corresponds to a decomposition $E=F_1\oplus F_2$ given by $F_1=\ker(A-\id_E)$ and $F_2=\ker(A+\id_E)$.
	The symmetric spaces
	$${\cal D} = \{(F_1,F_2)\in \Gr(E)^2\colon E=F_1\oplus F_2\}$$
	and $\Invol(\calL(E))$ are isomorphic.
\end{remark}
The following proposition shows that the tangent bundle of a symmetric space carries a natural symmetric space structure.
\begin{proposition}[{cf.\ \cite[Prop.~3.3]{Nee02}}] \label{prop:functorT}
	Let $(M,\mu)$ be a symmetric space and identify $T(M\times M)$ with $TM\times TM$. Then $(TM,T\mu)$ is a symmetric space. In each tangent space $T_xM$, $x\in M$, the product satisfies $v\cdot w = 2v-w$. For every morphism $f\colon (M_1,\mu_1)\rightarrow (M_2,\mu_2)$ of symmetric spaces, the tangent map $Tf$ is a morphism from $(TM_1,T\mu_1)$ into $(TM_2,T\mu_2)$. Thus $T$ is a covariant endofunctor of the category of symmetric spaces.
\end{proposition}
A smooth vector field $\xi\colon M \rightarrow TM$ is called a \emph{derivation} if it is a morphism of symmetric spaces. This can be rephrased by saying that $\xi\times \xi$ and $\xi$ are $\mu$-related vector fields. We denote the set of all such derivations by $\Der(M,\mu)$. By the naturality of the bracket product, it can be easily checked that $\Der(M,\mu)$ is a Lie subalgebra of the Lie algebra $\calV(M)$ of smooth vector fields (cf.\ \cite[I.5.7]{Ber08}).
%
%
%
\subsection{The Lie Triple System of a Pointed Symmetric Space}
\label{sec:lieTripleSystem}
	A \emph{Banach--Lie triple system}, simply called a \emph{Lie triple system}, is a Banach space $\mathfrak{m}$ with a continuous trilinear map 
	$[\cdot,\cdot,\cdot]\colon \mathfrak{m}^3 \rightarrow \mathfrak{m}$
	such that
	\begin{enumerate}
		\item[(1)] $[x,x,y]=0$,
		\item[(2)] $[x,y,z] + [y,z,x] + [z,x,y]=0$ \quad \mbox{and}
		\item[(3)] $[x,y,[u,v,w]]=[[x,y,u],v,w]+[u,[x,y,v],w]+[u,v,[x,y,w]]$
	\end{enumerate}
	for all $x,y,z,u,v,w \in \mathfrak{m}$.
	Given two Lie triple systems $\mathfrak{m}_1$ and $\mathfrak{m}_2$, a continuous linear map $A\colon \mathfrak{m}_1 \rightarrow \mathfrak{m}_2$ is called a \emph{morphism of Lie triple systems} if it satisfies
	$A[x,y,z]=[Ax,Ay,Az]$
	for all $x,y,z\in \mathfrak{m}_1$.

\begin{proposition} \label{prop:lieAlgebraLts}
	Let $\mathfrak{g}$ be a Banach--Lie algebra and $\sigma$ an involutive automorphism of $\mathfrak{g}$, i.e., $\sigma^2=\id_\mathfrak{g}$.
	Then $\mathfrak{g}$ can be written as the direct sum $\mathfrak{g}=\mathfrak{g}_+\oplus \mathfrak{g}_-$ of the Banach spaces 
	$\mathfrak{g}_+:=\ker(\sigma-\id)$
	and $\mathfrak{g}_-:=\ker(\sigma+\id)$.
	We have the rules
	$$[\mathfrak{g}_+,\mathfrak{g}_+]\subseteq \mathfrak{g}_+, \quad
		[\mathfrak{g}_+,\mathfrak{g}_-]\subseteq \mathfrak{g}_-, \quad\mbox{and}\quad
		[\mathfrak{g}_-,\mathfrak{g}_-]\subseteq \mathfrak{g}_+,$$
	and $\mathfrak{g}_-$ becomes a Banach--Lie triple system by
	$$[x,y,z]:=[[x,y],z].$$
\end{proposition}
\begin{proof}
	In \cite[p.~78]{Loo69}, O.~Loos states this in case of Lie triple systems and Lie algebras not equipped with a topological structure. In our case, the assertion remains true, as we simply have to check the continuity of $[\cdot,\cdot,\cdot]$, which immediately follows by the continuity of $[\cdot,\cdot]$ on the Banach--Lie algebra.
\end{proof}
The following proposition demonstrates how each tangent space of a symmetric space inherits a Lie triple structure by an isomorphism with a subspace of $\Der(M,\mu)$.
\begin{proposition}[{cf.\ \cite[I.5.9]{Ber08}}] \label{prop:killingSequence}
	Let $(M,\mu,b)$ be a pointed symmetric space. The symmetry $\mu_b$ around $b$ induces an involutive automorphism $(\mu_b)_\ast$ of the Lie algebra $\Der(M,\mu)$ by
	$$\begin{array}{cccc}
		(\mu_b)_\ast\colon & \Der(M,\mu) & \rightarrow & \Der(M,\mu) \\
		& \xi & \mapsto & T\mu_b \circ \xi \circ \mu_b.
	\end{array}$$
	Denoting the $(\pm 1)$-eigenspaces of $(\mu_b)_\ast$ by $\Der(M,\mu)_\pm$, we obtain a short exact sequence
	$$0 \rightarrow \Der(M,\mu)_+ \rightarrow \Der(M,\mu)=\Der(M,\mu)_+ \oplus \Der(M,\mu)_-
		\begin{temparraystretch}{0.0}
			\begin{array}{c} \rightarrow \\ \leftarrow \end{array}
		\end{temparraystretch}
	 	T_bM \rightarrow 0$$
	of vector spaces where the map $\Der(M,\mu)\rightarrow T_bM$ is the evaluation $\ev_b\colon \xi\mapsto\xi(b)$. It is split by the map
	$$\begin{array}{ccccc}
		T_bM & \stackrel{\cong}{\rightarrow} & \Der(M,\mu)_- &\subseteq& \Der(M,\mu) \\
		v & \mapsto & \xi_v
	\end{array}$$
	where $\xi_v(x):=\frac{1}{2} T\mu(v,T\mu(0_b,0_x))$.
	In particular, we have
	$$\Der(M,\mu)_+ = \{\xi\in \Der(M,\mu)\colon \xi(b)=0\} \quad\mbox{and}\quad
		\Der(M,\mu)_- = \{\xi_v\colon v\in T_bM\}.$$
\end{proposition}
\begin{theorem}[{cf.\ \cite[I.5.10]{Ber08}}]
	Let $(M,\mu,b)$ be a pointed symmetric space.
	By means of the linear isomorphism $\ev_b|_{\Der(M,\mu)_-}$ (cf.\ Proposition~\ref{prop:killingSequence}), the tangent space $T_bM$ becomes a Banach-Lie triple system. Given a morphism $f\colon (M_1,\mu_1,b_1) \rightarrow (M_2,\mu_2,b_2)$ of pointed symmetric spaces, the tangent map $T_{b_1}f$ is a morphism  between the respective Lie triple systems.
\end{theorem}
Equipped with this additional structure, $T_bM$ is called the \emph{Lie triple system of $(M,\mu,b)$}. We denote it by $\Lts(M,\mu,b)$ and in shorthand by $\mathfrak{m}$. The tangent map $T_{b_1}f$ of the morphism $f\colon (M_1,\mu_1,b_1) \rightarrow (M_2,\mu_2,b_2)$ is denoted by $\Lts(f)\colon \Lts(M_1,\mu_1,b_1) \rightarrow \Lts(M_2,\mu_2,b_2)$.
The assignment $\Lts$ is a functor from the category of pointed symmetric spaces to the category of Lie triple systems. This functor is called the \emph{Lie functor}.
%
%
%
\subsection{The Canonical Affine Connection of a Symmetric Space}
\label{sec:canonicalConnection}
For every symmetric space $(M,\mu)$, we have a canonical affine connection $B$ such that all symmetries $\mu_x$, $x\in M$, are affine automorphisms (cf.\ \cite[V.26]{Ber08}).
Given any chart \linebreak $\varphi\colon U \rightarrow V\subseteq E$ of $M$ and a restriction $\varphi_1\colon U_1 \rightarrow V_1\subseteq E$ of $\varphi$ such that $\mu(U_1\times U_1)\subseteq U$,\linebreak we denote the multiplication in these charts by $\mu^{\varphi_1}\colon V_1\times V_1\rightarrow V$. Then, the local representation $B^{\varphi_1}$ is given by
\begin{equation} \label{eqn:canonicalConnectionChart}
	B^{\varphi_1}_x(v,w) \ =\ \frac{1}{2}d^2\mu^{\varphi_1}(x,x)((0,v),(0,w)) \ =\ -\frac{1}{2}d^2\mu^{\varphi_1}(x,x)((v,0),(0,w)).
\end{equation}
In particular, this affine connection is torsionfree. (Thus it can be described by using a spray, cf.\ \cite[Th.~3.6]{Nee02}.)
\begin{theorem} \label{th:canonicalConnection}
	Let $(M,\mu)$ and $(M_2,\mu_2)$ be symmetric spaces and $\nabla$ and $\nabla\!_2$ the canonical affine connections on $TM$ and $TM_2$, respectively. Then the following assertions hold:
	\begin{enumerate}
		\item[\rm (1)] Let $f\colon M \rightarrow M_2$ be a smooth map. If it is a morphism of symmetric spaces then it is affine. If it is affine and $M$ is connected, then it is a morphism of symmetric spaces.
		\item[\rm (2)] $\Aut(M,\mu)= \Aut(M,\nabla)$, if $M$ is connected.
		\item[\rm (3)] $(M,\nabla)$ is geodesically complete.
		\item[\rm (4)] Let $\alpha\colon \RR\rightarrow M$ be a geodesic. For all $s,t\in\RR$, we have
		$$T_{\alpha(t+s)}\mu_{\alpha(t)}=-P_{t+s}^{t-s}(\alpha)\colon
			\quad T_{\alpha(t+s)}M \rightarrow  T_{\alpha(t-s)}M.$$
		\item[\rm (5)] Let $\alpha\colon \RR\rightarrow M$ be a geodesic and call the maps $\tau_{\alpha,s}:=\mu_{\alpha(\frac{1}{2}s)}\circ \mu_{\alpha(0)}$, $s\in \RR$, \emph{translations along $\alpha$}. Then these are automorphisms of $(M,\mu)$ with
		$$\tau_{\alpha,s}(\alpha(t))=\alpha(t+s) \quad\mbox{and}\quad
			T_{\alpha(t)}\tau_{\alpha,s}=P_{t}^{t+s}(\alpha)$$
		for all $s,t\in \RR$.
		\item[\rm (6)] Given a geodesic $\alpha\colon \RR\rightarrow M$, the vector field
		$$\xi_{\alpha}(x)
			:=\left.\frac{d}{dt}\right|_{t=0}\tau_{\alpha,t}(x)
			= \left.\frac{d}{dt}\right|_{t=0}
				\alpha(\textstyle\frac{1}{2}t)\cdot(\alpha(0)\cdot x)$$
		is an infinitesimal affine automorphism, and
		$\tau_{\alpha}\colon(t,x)\mapsto\tau_{\alpha,t}(x)$
		is its flow. The geodesic $\alpha$ is an integral curve of $\xi_{\alpha}$.
		\item[\rm (7)] Given a vector $v$ in $T_bM$, $b\in M$, let $\xi_v$ be the vector field from Proposition~\ref{prop:killingSequence} and $\alpha$ the maximal geodesic with $\alpha^\prime(0)=v$. Then we have $\xi_v=\xi_{\alpha}$.
	\end{enumerate}
\end{theorem}
\begin{proof}
	The assertions (2), (3) and (5) are due to \cite[Th.~3.6]{Nee02} in case of connected manifolds, but where the connectedness is not necessary for (3) and (5). The proof of (1) works in almost the same manner as the one of (2). Note that $\cite{Nee02}$ works with sprays. The assertions (4) - (6) follow from \cite[XIII, Prop.~5.3~-~5.5]{Lan01} whose proofs do also work in our context. Finally, we show (7) by an easy computation:
	\begin{eqnarray*}
		\xi_{\alpha}(x)
		&=& \textstyle\left.\frac{d}{dt}\right|_{t=0}
				\alpha(\frac{1}{2}t)\cdot(\alpha(0)\cdot x)
		\ = \ T\mu(\left.\frac{d}{dt}\right|_{t=0}\alpha(\frac{1}{2}t),
			0_{b \cdot x})
		\ = \ T\mu(\frac{1}{2}\alpha^\prime(0), T\mu(0_b, 0_x)) \\
		&=& \textstyle \frac{1}{2} T\mu(v,T\mu(0_b, 0_x))
		\ = \ \xi_v(x).
	\end{eqnarray*}
\end{proof}
\begin{remark}
	The preceding theorem shows that symmetric spaces satisfy the conditions for an affine locally symmetric space in the sense of S.~Helgason (cf.\ \cite[Ch.~IV, \S 1]{Hel01}). The symmetries coincide with the geodesic symmetries.
\end{remark}
\begin{theorem}[{cf.\ \cite[p.~137]{Ber08}}]
	Let $(M,\mu,b)$ be a pointed symmetric space with Lie triple system $\mathfrak{m}$ and curvature tensor $R$. Then we have
	$$[v,w,z] \ =\ -R_b(v,w,z)$$
	for all $v,w,z\in \mathfrak{m}$.
\end{theorem}
\begin{corollary}\label{cor:lieTipleMorphisms}
	Let $(M_1,\mu_1,b_1)$ and $(M_2,\mu_2,b_2)$ be pointed symmetric spaces, $R_1$ and $R_2$ their curvature tensors and $\mathfrak{m}_1$ and $\mathfrak{m}_2$ their Lie triple systems, respectively.
	Then a continuous linear map $A\colon \mathfrak{m}_1 \rightarrow \mathfrak{m}_2$ is a morphism of Lie triple systems if and only if it intertwines the curvature tensors in $b_1$ and $b_2$ in the sense that
	$$A\circ (R_1)_{b_1} \ =\ (R_2)_{b_2}\circ A^3.$$
\end{corollary}
\begin{proposition}\label{prop:parallelR}
	Given a symmetric space $(M,\mu)$, its curvature tensor $R$ is parallel on $M$, i.e., $\nabla\!_uR=0$ for all vectors $u\in TM$.
\end{proposition}
\begin{proof}
	Cf.\ the proof of \cite[Prop.~XIII.6.2]{Lan01}.
\end{proof}
\begin{proposition}[{cf.\ \cite[p.~84]{Loo69}}] \label{prop:Der=Kill}
	Let $(M,\mu)$ be a connected symmetric space. A smooth vector field $\xi$ on $M$ is a derivation if and only if it is an infinitesimal affine automorphism with respect to the canonical affine connection on $TM$, i.e., we have $\Der(M,\mu) = \Kill(M,\nabla)$.
\end{proposition}
%
\begin{proof}
	$\Der(M,\mu)\subseteq \Kill(M,\nabla)$: Given a derivation $\xi\in \Der(M,\mu)$, we shall check that for each point $b$ of $M$, there is a chart $\varphi_1\colon U_1\rightarrow V_1\subseteq E$ at $b$ such that the equation (\ref{eqn:infAffAutoChart}) of Section~\ref{sec:affAndInfAffAuto} is satisfied.
	For this, we consider some chart $\varphi\colon U\rightarrow V\subseteq E$ at $b$ and let $\varphi_1\colon U_1\rightarrow V_1\subseteq E$ be a restriction of $\varphi$ such that $b\in U_1$ and $\mu(U_1\times U_1)\subseteq U$. We denote the multiplication in these charts by $\mu^{\varphi_1}\colon V_1\times V_1\rightarrow V$. The vector field $\xi$ being a morphism, we have
	$\xi\circ\mu=T\mu\circ(\xi\times\xi)$,
	which locally means
	$$\xi^\varphi(\mu^{\varphi_1}(x,y))
		=d\mu^{\varphi_1}(x,y)(\xi^\varphi(x),\xi^\varphi(y))$$
	for all $x,y\in V_1$.
	By taking the derivative with respect to $x$ in any direction $v\in E$, we obtain
	\begin{eqnarray*}
		d\xi^\varphi(\mu^{\varphi_1}
			(x,y))(d\mu^{\varphi_1}(x,y)(v,0))
		&=& d^2\mu^{\varphi_1}(x,y)
			\big((\xi^\varphi(x),\xi^\varphi(y)),(v,0)\big) \\
		&&\quad {}+ d\mu^{\varphi_1}(x,y)(d\xi^\varphi(x)(v),0).
	\end{eqnarray*}
	We then take the derivative with respect to $y$ in any direction $w\in E$ and get
	\begin{eqnarray*}
		\lefteqn{d^2\xi^\varphi(\mu^{\varphi_1}
			(x,y))\big(d\mu^{\varphi_1}(x,y)(v,0),
				d\mu^{\varphi_1}(x,y)(0,w)\big)} \\
		&& \quad {}+ d\xi^\varphi(\mu^{\varphi_1}(x,y))
				\big(d^2\mu^{\varphi_1}(x,y)((v,0),(0,w))\big) \\
		&\quad =& d^3\mu^{\varphi_1}(x,y)
			\big((\xi^\varphi(x),\xi^\varphi(y)),(v,0),(0,w)\big)
			+ d^2\mu^{\varphi_1}(x,y)
			\big((0,d\xi^\varphi(y)(w)),(v,0)\big) \\
		&&\quad {}+ d^2\mu^{\varphi_1}(x,y)
			\big((d\xi^\varphi(x)(v),0),(0,w)\big).
	\end{eqnarray*}
	By putting $y:=x$ and applying (\ref{eqn:canonicalConnectionChart}) and Proposition~\ref{prop:functorT}, we obtain
	\begin{eqnarray*}
		\lefteqn{d^2\xi^\varphi(x)(2v,-w)
			+ d\xi^\varphi(x)
				(-2B^{\varphi_1}_x(v,w))} \\
		&\quad =& -2dB^{\varphi_1}(x)(\xi^\varphi(x))(v,w)
			- 2B^{\varphi_1}_x(v,d\xi^\varphi(x)(w))
			- 2B^{\varphi_1}_x(d\xi^\varphi(x)(v),w)
	\end{eqnarray*}
	for all $x\in V_1$ and $v,w\in E$, which entails the required equation when multiplying the two sides by $(-\frac{1}{2})$.
		
	$\Kill(M,\nabla)\subseteq \Der(M,\mu)$: Given an infinitesimal automorphism $\xi$, we know from Theorem~\ref{th:canonicalConnection}(3) and Theorem~\ref{th:Kill(M,nabla)Complete} that it is complete.
	Being affine automorphisms, the flow maps $\flow^\xi_t\colon M\rightarrow M$ are also automorphisms of symmetric spaces, by Theorem~\ref{th:canonicalConnection}(2). Therefore we have
	\begin{eqnarray*}
		\xi(x\cdot y)
		&=& \textstyle \left.\frac{d}{dt}\right|_{t=0}\flow^\xi_t(x\cdot y)
		\ = \ \left.\frac{d}{dt}\right|_{t=0}\flow^\xi_t(x)\cdot\flow^\xi_t(y)
		\ = \ T\mu(\left.\frac{d}{dt}\right|_{t=0}\flow^\xi_t(x),
			\left.\frac{d}{dt}\right|_{t=0}\flow^\xi_t(y)) \\
		&=& T\mu(\xi(x),\xi(y))
		\ = \ \xi(x)\cdot \xi(y),
	\end{eqnarray*}
	hence $\xi$ is a derivation.
\end{proof}
\begin{corollary}
	Let $(M,\mu)$ be a connected symmetric space. The Lie algebra $\Der(M,\mu)$ of derivations can be turned into a Banach--Lie algebra according to Proposition~\ref{prop:isomorphismOfKill}: Its Banach space structure is uniquely determined by the requirement that for each $p\in\Fr(M)$, the map
	$$\Der(M,\mu)\rightarrow T_p(\Fr(M)),\ \xi\mapsto \left.\frac{d}{dt}\right|_{t=0}\Fr(\flow^{\xi}_t)(p)$$
	is a closed embedding.
\end{corollary}
\begin{corollary}
	Let $(M,\mu)$ be a symmetric space. Every derivation $\xi\in\Der(M,\mu)$ is complete.
\end{corollary}
\begin{proof}
	We can without loss of generality assume that $(M,\mu)$ is connected, as the matter of local flows takes place in connected components. We then know from Theorem~\ref{th:canonicalConnection}(3) and Theorem~\ref{th:Kill(M,nabla)Complete} that $\xi$ is complete.
\end{proof}
%
%
%
%
\subsection{Integration of Morphisms of Lie Triple Systems}
\label{sec:integrationLTS}
In the light of Theorem~\ref{th:canonicalConnection}(1), it is clear that two morphisms $f,g\colon (M_1,\mu_1,b_1)\rightarrow  (M_2,\mu_2,b_2)$ of pointed symmetric spaces that satisfy $\Lts(f)=\Lts(g)$ are equal if $M_1$ is connected. Considering morphisms of Lie triple systems, the following theorem deals with the existence of integrals.

\begin{theorem}[Integrability Theorem]
	Let $(M_1,\mu_1,b_1)$ and $(M_2,\mu_2,b_2)$ be pointed symmetric spaces with Lie triple systems $\mathfrak{m}_1$ and $\mathfrak{m}_2$, respectively, and let $A\colon \mathfrak{m}_1 \rightarrow \mathfrak{m}_2$ be a morphism of Lie triple systems.
	If $M_1$ is 1-connected, i.e., connected and simply connected, then there exists a unique morphism $f\colon M_1\rightarrow M_2$ of pointed symmetric spaces that satisfies $\Lts(f)=A$.
\end{theorem}
\begin{proof}
	We equip the tangent bundles with their canonical affine connections and apply Theorem~\ref{th:globalIntegrability}, Theorem~\ref{th:canonicalConnection}, Corollary~\ref{cor:lieTipleMorphisms} and Proposition~\ref{prop:parallelR}.
\end{proof}
%
%
%
%
\subsection{The Automorphism Group of a Symmetric Space}
\label{sec:autoOfSymSpace}
\begin{lemma}\label{lem:transitiveAction}
	Let $(M,\mu)$ be a connected symmetric space. The automorphism group\linebreak $\Aut(M,\mu)$ acts transitively on $M$.
\end{lemma}
\begin{proof}
	Due to Theorem~\ref{th:canonicalConnection}, there are translations along geodesics. As $M$ is geodesically connected (cf.\ Lemma~\ref{lem:geodesicConnection}), each point of $M$ can be mapped to any other one by composing translations.
\end{proof}
\begin{theorem}\label{th:autoGroupOfSymSpace}
	Let $(M,\mu)$ be a connected symmetric space.
	The automorphism group $\Aut(M,\mu)$ can be turned into a Banach--Lie group such that
	$$\exp\colon \Der(M,\mu)\rightarrow \Aut(M,\mu),\ \xi \mapsto \flow^{-\xi}_1$$
	is its exponential map. The natural map $\tau\colon \Aut(M,\mu)\times M \rightarrow M$ is a transitive smooth action whose derived action is the inclusion $\Der(M,\mu)\hookrightarrow\calV(M)$, i.e., $-T\tau(\id_M,x)(\xi,0)=\xi(x)$.
\end{theorem}
\begin{proof}
	This theorem is a direct consequence of Theorem~\ref{th:autoGroupOfAffineManifold}, Theorem~\ref{th:canonicalConnection}, Proposition~\ref{prop:Der=Kill} and Lemma~\ref{lem:transitiveAction}.
\end{proof}
\begin{proposition}\label{prop:killingSequence2}
	Let $(M,\mu,b)$ be a pointed connected symmetric space. The short exact sequence
	$$0 \rightarrow \Der(M,\mu)_+ \rightarrow \Der(M,\mu)=\Der(M,\mu)_+ \oplus \Der(M,\mu)_-
		\begin{temparraystretch}{0.0}
			\begin{array}{c} \rightarrow \\ \leftarrow \end{array}
		\end{temparraystretch}
	 	T_bM \rightarrow 0$$
	of vector spaces (cf.\ Proposition~\ref{prop:killingSequence}) is actually one of Banach spaces.
\end{proposition}
\begin{proof}
	The evaluation map $\ev_b\colon \Der(M,\mu)\rightarrow T_bM$ is continuous by Theorem~\ref{th:autoGroupOfSymSpace} and hence, $\Der(M,\mu)_+=\ker(\ev_b)$ is a closed subspace of $\Der(M,\mu)$. To see that $\Der(M,\mu)_-=\ker((\mu_b)_\ast+\id_{\Der(M,\mu)})$ is closed, too, we shall check that
	$$(\mu_b)_\ast\colon \Der(M,\mu)\rightarrow \Der(M,\mu),\ \xi\mapsto T\mu_b\circ\xi\circ\mu_b$$
	is continuous and shall do this by proving $L(c_{\mu_b})=(\mu_b)_\ast$, where
	$$c_{\mu_b}\colon \Aut(M,\mu)\rightarrow \Aut(M,\mu),\ g\mapsto\mu_b\circ g \circ \mu_b$$
	denotes the conjugation map.
	We have to verify that $\mu_b\circ\flow^{-\xi}_t\circ\mu_b=\flow^{-(\mu_b)_\ast\xi}_t$ for all $t\in\RR$.
	This is true by the general fact that the flow of the push-forward of a vector field is the push-forward of its flow (cf.\ \cite[Prop.~4.2.4]{AMR88}).
\end{proof}
\begin{theorem}[Homogeneity of connected symmetric spaces]
	Let $(M,\mu)$ be a connected symmetric space with automorphism group $G:=\Aut(M,\mu)$. Given $b\in M$, the stabilizer subgroup $G_b\leq G$ is an open subgroup of the group $G^{c_{\mu_b}}$ of fixed points in $G$ for the involution $c_{\mu_b}$ on $G$ given by $c_{\mu_b}(g):=\mu_b\circ g\circ \mu_b$. The symmetric space $G/G_b$ (cf.\ Example~\ref{ex:homogeneousSpaces}) is ismorphic to $M$ via the isomorphism $\Phi\colon G/G_b\rightarrow M$ given by $\Phi(gG_b):=g(b)$.	 
\end{theorem}
\begin{proof}
	Given $g\in G$, we have $g\circ\mu_b = \mu_{g(b)}\circ g$, so that $g\in G^{c_{\mu_b}}$ if $g(b)=b$. Hence, the stabilizer $G_b$ is a subgroup of $G^{c_{\mu_b}}$.
	
	To see that $G_b$ is open in $G^{c_{\mu_b}}$, we shall prove that $G_b=\tau_b^{-1}(U) \cap G^{c_{\mu_b}}$ for the orbit map $\tau_b\colon\Aut(M,\mu)\rightarrow M, \ g\mapsto g(b)$, where $U$ is an open neighborhood of $b$ such that $b$ is an isolated fixed point of $(\mu_b)|_U$. Given $g\in G^{c_{\mu_b}}$, the point $\tau_b(g)=g(b)$ is a fixed point of $\mu_b$, so that $\tau_b(g)\in U$ if and only if $g(b)=b$. This shows $G_b=\tau_b^{-1}(U) \cap G^{c_{\mu_b}}$.
	
	It is clear that the map $\Phi$ is well-defined and injective. Its surjectivity follows by the transitivity of the natural action of $G$ on $M$.
	
	Due to $\Phi\circ q = \tau_b$ with the quotient map $q\colon G\rightarrow G/G_b$, the map $\Phi$ is automatically smooth, as $q$ is a submersion. To see the smoothness of the inverse map $\Phi^{-1}$, we shall show that $\tau_b$ is a submersion. By Theorem~\ref{th:autoGroupOfSymSpace}, we have $T_{\id_M}\tau_b=-\ev_b$ with the evaluation map $\ev_b\colon \Der(M,\mu)\rightarrow T_bM$, so that $\tau_b$ is submersive at $\id_M$, the map $\ev_b$ being surjective with splitting kernel (cf.\ Proposition~\ref{prop:killingSequence2}). For every $g\in G$, we have
	$\tau_b = g\circ \tau_b \circ \lambda_{g^{-1}},$
	where $\lambda_{g^{-1}}$ is the left multiplication with $g^{-1}$ in $G$, whence $\tau_b$ is submersive at $g$, too. Therefore, $\tau_b$ is a submersion.
	
	Finally, we check that $\Phi$ is a homomorphism:
	\begin{eqnarray*}
		\Phi(gG_b \cdot hG_b) & =& \Phi(gc_{\mu_b}(g)^{-1}c_{\mu_b}(h)G_b) \ =\ (g \mu_bg^{-1}\mu_b \mu_bh\mu_b)(b) \ =\ (g\mu_bg^{-1})(h(b)) \\
		&=& g\big(b\cdot g^{-1}(h(b))\big) \ =\ g(b) \cdot g\big(g^{-1}(h(b))\big) \ =\ \Phi(gG_b)\cdot \Phi(hG_b).
	\end{eqnarray*}
\end{proof}
\begin{acknowledgements}
	I am grateful to Karl-Hermann Neeb for his helpful communications and proof reading during my research towards this article. This work was supported from the Technical University of Darmstadt and from the Studienstiftung des deutschen Volkes.
\end{acknowledgements}

\bibliography{paperB}
\bibliographystyle{amsalpha}
\end{document}